\documentclass[a4paper, 10pt]{article}
\usepackage[latin1]{inputenc}
\usepackage[english]{babel}

\usepackage{amsmath}
\usepackage{amsthm}
\usepackage{amssymb}
\usepackage{color}
\usepackage{bbm}
\usepackage{todonotes}
\usepackage{multirow}
\usepackage{graphicx}

\theoremstyle{plain}
\newtheorem{theorem}{Theorem}%[section]
\newtheorem{proposition}[theorem]{Proposition}

\theoremstyle{remark}
\newtheorem{remark}[theorem]{Remark}

\newcommand{\R}{\mathbb{R}}

\newcommand{\mD}{\mathcal{D}}

\newcommand{\mR}{\mathcal{R}}
\newcommand{\mT}{\mathcal{T}}
\newcommand{\mV}{\mathcal{V}}

\newcommand{\dd}{\, \text{d}}
\newcommand{\matlab}{MATLAB\textsuperscript{\textregistered}}
\newcommand{\intel}{Intel\textsuperscript{\textregistered}}

\usepackage{array}
\newcolumntype{R}[1]{>{\raggedleft\let\newline\\\arraybackslash\hspace{0pt}}m{#1}}

\usepackage{tikz}
\usepackage{pgfplots}
\usepackage{subcaption}
\usepackage{caption}

\usepackage{hyperref}
\usepackage{vmargin}
\setmarginsrb{3cm}{2.5cm}{3cm}{1.5cm}{0cm}{0cm}{0cm}{1.5cm}

\usetikzlibrary{external}
{FP_numerics}

\begin{document}

\title{Numerical Study of Polynomial Feedback Laws
for a Bilinear Control Problem}
\author{Tobias Breiten\footnote{Institute of Mathematics and Scientific Computing, University of Graz, Heinrichstrasse 36, 8010 Graz, Austria. E-mail: tobias.breiten@uni-graz.at} \quad
Karl Kunisch\footnote{Institute of Mathematics and Scientific Computing, University of Graz, Heinrichstrasse 36, 8010 Graz, Austria and RICAM Institute, Austrian Academy of Sciences, Altenbergerstrasse 69, 4040 Linz, Austria. E-mail: karl.kunisch@uni-graz.at} \quad
Laurent Pfeiffer\footnote{Institute of Mathematics and Scientific Computing, University of Graz, Heinrichstrasse 36, 8010 Graz, Austria. E-mail: laurent.pfeiffer@uni-graz.at}
}
\maketitle

\begin{abstract}
An infinite-dimensional bilinear optimal control problem with infinite-time horizon is considered.
The associated value function can be expanded in a Taylor series around the equilibrium, the Taylor series involving multilinear forms which are uniquely characterized by generalized Lyapunov equations.
A numerical method for solving these equations is proposed. It is based on a generalization of the balanced truncation model reduction method and some techniques of tensor calculus, in order to attenuate the curse of dimensionality.
Polynomial feedback laws are derived from the Taylor expansion and are numerically investigated for a control problem of the Fokker-Planck equation. Their efficiency is demonstrated for initial values which are sufficiently close to the equilibrium.
\end{abstract}



\section{Introduction}

In this article, we consider the following bilinear optimal control problem:
\begin{equation}\label{eq:setup}
\begin{aligned}
& \inf_{u \in L^2(0,\infty;\mathbb R^m)} \mathcal{J}(u,y_0) := \frac{1}{2}
\int_0^\infty \| y(t) \|_Y^2 \dd
t + \frac{\beta}{2} \int_0^\infty \|u(t)\|_{\mathbb{R}^m}^2 \dd t, \\[1.6ex]
& \qquad \text{where: }
  \left\{    \begin{array} {l} \displaystyle
\dot{y}(t)= Ay(t)+\sum_{j=1}^m(N_jy(t)+B_j)u_j(t), \quad
\text{for
} t>0 \\
y(0)=y_0 \in Y.
\end{array} \right.
\end{aligned}
\end{equation}
Here, $V \subset Y \subset V^*$ is a Gelfand triple of real Hilbert spaces and
$A\colon \mD(A)\subset Y \rightarrow Y$ is the infinitesimal generator of an
analytic $C_0$-semigroup $e^{At}$ on $Y$. The precise conditions on
$B_j$ and
$N_j$ are given further below. The value function, denoted by $\mathcal{V}$, associates with any initial condition $y_0$ the value of problem \eqref{eq:setup}.

In our previous work \cite{BreKP17b}, we analysed polynomial feedback laws of the form
\begin{equation*}
\big( \mathbf{u}_p(y) \big)_j = -\frac{1}{\beta} D\mathcal{V}_p(y)(N_j y + B_j), \quad \forall j=1,...,m
\end{equation*}
resulting from a Taylor expansion $\mathcal{V}_p$ of the value function $\mathcal{V}$.
The Taylor expansion is of the following form:
\begin{align}\label{eq:val_fun_pol_app}
\mathcal{V}(y) \approx \mathcal{V}_p(y)=\sum_{k=2}^p \frac{1}{k!}
\mathcal{T}_k(y,\dots,y),
\end{align}
where $\mathcal{T}_k\colon Y^k \rightarrow \R$ denotes a bounded multilinear form of order $k$.
The multilinear form $\mathcal{T}_2$ is determined by solving an algebraic Riccati equation. For $k \geq 3$, the multilinear form $\mathcal{T}_k$ is characterized by a generalized operator Lyapunov equation of the
form
\begin{align}\label{eq:operator_eq}
\sum_{i=1}^k \mathcal{T}_k(z_1,\dots,z_{i-1},A_\Pi z_i, z_{i+1},\dots,
z_k) =
\mathcal{R}_k(z_1,\dots,z_k), \ \ z_1,\dots,z_k \in \mathcal{D}(A).
\end{align}
In this equation, the operator $A_\Pi$ is associated with the
linearized closed-loop system obtained from a Riccati-based stabilization
approach and the r.h.s.\@ is computed by induction.

In the present contribution, we provide a detailed description
of the numerical implementation of
the feedback laws $\mathbf{u}_p$ and investigate their behavior
in practice.
A version of the Lyapunov equations \eqref{eq:operator_eq} in a finite-dimensional space is obtained by discretizing the state equation with a finite-difference scheme, which preserves the bilinear structure of the system. The numerical realisation of the discretized Lyapunov is not straightforward, because of the curse of dimensionality: The size of the linear system to be solved increases exponentially with the dimension of the domain $\Omega$ and with the degree $p$ of the Taylor expansion. We therefore propose to use a generalization of the balanced truncation model reduction method \cite{BenD11} to reduce the size of the dimension of the state equation and the Lyapunov equations. We also use a technique of \cite{Gra04} for solving the discretized and reduced Lyapunov equations.

The method is tested on an optimal control problem of the Fokker-Planck
equation, in dimension 1 and 2. The impact of model reduction on the
corresponding feedback laws is analysed. The efficiency of the feedback laws is
also analysed, in particular, we investigate how much improvement can be
obtained when using high-order feedback laws rather than Riccati-based feedback
laws.
At a theoretical level, the method is of local nature: The well-posedness of the closed-loop system associated with the feedback law $\mathbf{u}_p$ is only guaranteed for initial conditions close to the equilibrium \cite{BreKP17b}. We therefore investigate the influence of the distance of the initial condition to the origin on the success of the method. The influence of the cost parameter $\beta$ is also investigated.

The concept of series expansion of the value function has inspired
researchers in optimal control for a long time. We refer to \cite{KreAH13} for a very useful survey on
this topic. Numerical tests are mostly carried out in the context of
ordinary differential equations for systems of significantly smaller
order than those which typically arise from discretized infinite-dimensional systems, which are in the focus of the present work. Let us mention that an
interesting and natural extension of the expansion method consists in
computing the expansion in an array of reference points. This  relates
to the concept of the patchy technique \cite{AguK14,AncB99}. Many
additional references have been gathered in \cite{BreKP17b}.

Important efforts have been made recently to develop some new methods for the
feedback control of partial differential equations. The present paper
contributes to this developing field. For quadratic cost functionals and in the
absence of additional constraints, the most noted and investigated technique
consists in applying linear quadratic regulator theory after linearization of
the state equation around a steady state, see for instance \cite{B09,Ray06}.
Most other techniques, and especially those involving the value function and the
Hamilton-Jacobi-Bellmann equation, rely on system reduction. In
\cite{AF13,KVX04}, dimension reduction of the system is based on proper
orthogonal decomposition, while the results in \cite{KalK17} are based the
high-order approximation properties of spectral methods. In \cite{ItoR98}, the
reduced basis method was used for open-loop control of fluid flow.

The structure of the article is as follows. In Section
\ref{sec:recap_MIMO}, we
briefly recall the main theoretical results obtained in \cite{BreKP17b}
and generalize them to the case of multiple inputs.
In Section \ref{sec:FokkerPlanck}, we describe the bilinear control problem of
the Fokker-Planck equation used for the numerical results.
We provide in Section \ref{sec:alg_app} a detailed description of our numerical
approach for discretizing, reducing, and solving the Lyapunov equations.
Numerical results are reported in Section \ref{sec:num_res}.

\section{Summary of the theoretical results}
\label{sec:recap_MIMO}

In this section, we recall the main theoretical results proved in \cite{BreKP17b}.
In that previous paper, we worked with scalar-valued controls, while we now
consider the multi-input control case $u(t) \in \R^m$.
The extension of our results to the case $m>1$ is however straightforward.
Throughout this section, we assume that the following assumptions
are satisfied.

\begin{itemize}
 \item[(A1)] The operator $A$ can be associated with a bounded $V$-$Y$ bilinear
form $a\colon V\times V \to \mathbb R$ such that there exist $\lambda \in \mathbb
R$ and $\delta > 0$ satisfying
\begin{align*}
  a(v,v) \ge \delta \| v \| _V^2 - \lambda \| v\|_Y^2  \quad \text{for all } v \in V.
\end{align*}

\item[(A2)] For all $j=1,\dots,m$, $N_j \in \mathcal{L}(V,Y) \cap \mathcal{L}(\mathcal{D}(A), V)$ and
$N_j^* \in \mathcal{L}(V,Y)$.

\item[(A3)] For $A_0:=A-\gamma I,$ with $\gamma >0$ sufficiently large and the
real interpolation space with indices 2 and $\frac{1}{2}$, see
\cite[Proposition 6.1, Part II, Chapter 1]{Benetal07}, it holds that
\begin{align*}
  [\mathcal{D}(-A_0),Y]_{\frac{1}{2}} = [\mathcal{D}(-A_0^*),Y]_{\frac{1}{2}}
= V.
\end{align*}

\item[(A4)] There exist operators $F_1,\dots,F_m \in \mathcal{L}(Y,\mathbb R)$
such that the semigroup $e^{(A+\sum_{j=1}^m B_j F_j) t}$ is exponentially stable on $Y.$

\end{itemize}

In the approach developed in \cite{BreKP17b}, we first characterize the
multilinear forms $\mathcal{T}_2$,$\mathcal{T}_3,\dots$ of the Taylor
approximation. The equations satisfied by the multilinear forms are obtained by
successive differentiation of the following Hamilton-Jacobi-Bellman equation.

\begin{proposition}[Proposition 9, \cite{BreKP17b}]\label{prop:HJB}
Assume that there exists an open neighborhood $Y_0$ of the origin in $Y$ which is such that the two following
statements hold:
\begin{enumerate}
\item For all $y_0 \in Y_0$, problem \eqref{eq:setup} possesses a solution $u$
which is right-continuous at time 0.
\item The value function is continuously differentiable on $Y_0$.
\end{enumerate}
Then, for all $y_0 \in \mathcal{D}(A) \cap Y_0$,
\begin{equation}\label{eq:HJB}
D\mV(y) Ay_0 + \frac{1}{2} \|y_0 \|_Y^2 - \frac{1}{2\beta} \sum_{j=1}^m \big(
D\mV(y_0) (N_jy_0+B_j) \big)^2=0.
\end{equation}
Moreover, for all solutions $u$ to problem \eqref{eq:setup} with initial condition $y_0$, if $u$ is right-continuous at 0, then
\begin{equation} \label{eqFeedbackGeneral}
\big(u(0^+) \big)_j= - \frac{1}{\beta} D\mathcal{V}(y_0)(N_j y + B_j), \quad \forall j=1,\dots,m
\end{equation}
\end{proposition}

The first nontrivial term $\mathcal{T}_2$ of the Taylor expansion
is determined by the unique nonnegative self-adjoint operator satisfying the
following
algebraic operator Riccati equation:
\begin{equation} \label{eq:algebraicRiccati}
\langle A^* \Pi z_1,z_2 \rangle + \langle \Pi A z_1,z_2 \rangle + \langle
z_1,z_2\rangle - \frac{1}{\beta} \sum_{j=1}^m (B_j^* \Pi z_1)(B_j^* \Pi z_2)=
0 \ \ \text{for all } z_1,z_2 \in \mD(A).
\end{equation}
It is obtained by differentiating twice equation \eqref{eq:HJB}.
It is well-known that the linearized closed-loop operator
\begin{equation}
A_\Pi : = A - \frac{1}{\beta} \sum_{j=1}^m B_j B_j^* \Pi
\end{equation}
generates an exponentially stable semigroup on $Y$, thanks to assumption (A4).

Further differentiation of the HJB equation allows to characterize
the multilinear forms $\mathcal{T}_3$, $\mathcal{T}_4$... as solutions to
generalized Lyapunov equations, whose right-hand sides are defined recursively.
The precise structure of these equations is given in Theorem \ref{thm:mult_lin_form} below.
In the definition of the right-hand sides of the generalized Lyapunov equations, we make use
of a specific symmetrization technique, that we define now.
For $i$ and $j \in \mathbb{N}$, consider the following set of permutations:
\begin{equation*}
S_{i,j}= \big\{ \sigma_{i+j} \, | \, \sigma(1) < ... < \sigma(i) \text{ and } \sigma(i+1) < ... < \sigma(i+j) \big\},
\end{equation*}
where $S_{i,j}$ is the set of permutations of $\{ 1,...,i+j \}$.
Let $\mathcal{T}$ be a multilinear form of order $i+j$. We denote by $\text{Sym}_{i,j}(\mathcal{T})$ the multilinear form defined by
\begin{equation*}
\text{Sym}_{i,j}(\mathcal{T})(z_1,...,z_{i+j})
= \binom{i+j}{i}^{-1} \Big[ \sum_{\sigma \in S_{i,j}}
\mathcal{T}(z_{\sigma(1)},...,z_{\sigma(i+j)} ) \Big], \quad \forall (z_1,...,z_{i+1}) \in Y^{i+j}.
\end{equation*}

\begin{theorem}[Theorem 15, \cite{BreKP17b}] \label{thm:mult_lin_form}
There exists a unique sequence of bounded symmetric multilinear forms $(\mT_k)_{k \geq
2,}$ with $\mT_k\colon Y^k \rightarrow \R$ and a unique sequence of bounded multilinear forms $(\mR_{j,k})_{k \geq
3},j=1,\dots,m$
with $\mR_{j,k}\colon \mD(A)^k \rightarrow \R$ such that for all $(z_1,z_2) \in Y^2$,
\begin{equation} \label{eqLyapounov3d}
\mathcal{T}_2(z_1,z_2):= (z_1,\Pi z_2)
\end{equation}
and such that for  all $k \geq 3$, for all $(z_1,...,z_k) \in \mathcal{D}(A)^k$,
\begin{subequations} \label{eqLyapounov3}
\begin{equation}
\sum_{i=1}^k \mathcal{T}_k (z_1,...,z_{i-1},A_\Pi z_i, z_{i+1},...,z_k)=
\frac{1}{2\beta} \sum_{j=1}^m\mathcal{R}_{j,k}(z_1,...,z_k), \label{eqLyapounov3a}
\end{equation}
where
\begin{align}
\mathcal{R}_{j,k}= \ & 2k(k-1) \text{\emph{Sym}}_{1,k-1}\big( \mathcal{C}_{j,1}
\otimes
\mathcal{G}_{j,k-1} \big) \notag \\
& \qquad + \sum_{i=2}^{k-2} \binom{k}{i}
\text{\emph{Sym}}_{i,k-i} \big( (\mathcal{C}_{j,i} + i \mathcal{G}_{j,i})
\otimes
(\mathcal{C}_{j,k-i} + (k-i) \mathcal{G}_{j,k-i}) \big), \label{eqLyapounov3b}
\end{align}
and where
\begin{equation} \label{eqLyapounov3c}
\begin{cases} \begin{array}{rl}
\mathcal{C}_{j,i}(z_1,...,z_i) = & \mathcal{T}_{i+1}(B_j,z_1,...,z_i),
\quad \text{for $i=1,...,k-2$,}\\
\mathcal{G}_{j,i}(z_1,...,z_i) = &
\frac{1}{i} \Big[
\sum_{\ell=1}^i \mathcal{T}_i(z_1,...,z_{\ell-1},N_jz_\ell,z_{\ell+1},...,z_i) \Big],
\quad \text{for $i=1,...,k-1$}.
\end{array} \end{cases}
\end{equation}
\end{subequations}
\end{theorem}

For all $p \geq 2$, we define the polynomial approximation $\mathcal{V}_p$ as follows:
\begin{equation}\label{eq:kk20}
  \begin{aligned}
    \mathcal{V}_p & \colon Y \to \R,\quad
    \mathcal{V}_p(y) = \sum_{k=2}^p \frac{1}{k!} \mathcal{T}_k(y,\dots,y),
  \end{aligned}
  \end{equation}
where the sequence $(\mathcal{T}_k)_{k \geq 2}$ is given by Theorem \eqref{thm:mult_lin_form}.
We deduce from $\mathcal{V}_p$ the polynomial feedback law $\mathbf{u}_p \colon y \in V \rightarrow \R^m$, defined by
\begin{align}
\big(\mathbf{u}_p(y)\big)_j = \ & -
\frac{1}{\beta} D \mathcal{V}_p(y)(N_j y+B_j) \notag \\
= \ & - \frac{1}{\beta} \Big( \sum_{k=2}^p \frac{1}{(k-1)!} \mathcal{T}_k(N_j y +
B_j,y,\dots,y ) \Big), \quad \forall j=1,\dots,m. \label{eq:feedback_law}
\end{align}
Its form is suggested by \eqref{eqFeedbackGeneral} and \eqref{eq:kk20}.
A justification of the differentiability of $\mathcal{V}_p$ and a formula for its derivative, used in the above expression, can be found in \cite[Lemma 7]{BreKP17b}.
We consider now the closed-loop system associated with the feedback law $\mathbf{u}_p$:
\begin{equation}\label{eq:cl_poly}
\dot{y}(t) = Ay(t) + \sum_{j=1}^m (N_j y(t)+ B_j) \big( \mathbf{u}_p(y(t))\big)_j, \quad
y(0) = y_0.
\end{equation}
For a given initial condition $y_0$, its solution is denoted by $S(\mathbf{u}_p,y_0)$.
We also denote by $\mathbf{U}_p(y_0)$ the open-loop control defined by
\begin{equation} \label{eqDefClosedLoop}
\mathbf{U}_p(y_0;t)= \mathbf{u}_p(S(\mathbf{u}_p,y_0;t)), \quad \text{for a.e.\@ $t \geq 0$.}
\end{equation}
%We call it control generated by $\mathbf{u}_p$ and $y_0$.
The following theorem states that for $\| y_0 \|_Y$ small enough,
the closed-loop system \eqref{eq:cl_poly} has a unique solution
and generates an open-loop control in $L^2(0,\infty;\R^m)$.
The solution to the closed-loop system is obtained in the space:
\begin{equation*}
W_\infty:= \Big\{ y \in L^2(0,\infty;V) \,|\, \dot{y} \in L^2(0,\infty;V^*)
\Big\}.
\end{equation*}

\begin{theorem}[Theorem 21 and Corollary 22, \cite{BreKP17b}] \label{thm:well_posed_feedback}
There exist two constants $\delta_0>0$ and $C>0$
such that for all $y_0$ with $\| y_0 \|_Y \leq \delta_0$, the closed-loop system
\eqref{eq:cl_poly} admits a unique solution $S(\mathbf{u}_p,y_0) \in W_\infty$ satisfying
\begin{equation}
\| S(\mathbf{u}_p,y_0) \| _{W_\infty} \leq C \| y_0 \|_Y,
\end{equation}
moreover, $\mathbf{U}_p(y_0) \in L^2(0,\infty;\R^m)$.
\end{theorem}

Finally, the following theorem states that
$\mathcal{V}_p$ is an approximation of $\mathcal{V}$ of order $p+1$,
in the neighborhood of 0 and gives an error estimate on the efficiency of
the open-loop control generated by $\mathbf{u}_p$.

\begin{theorem}[Proposition 2, Theorem 30, and Theorem 32, \cite{BreKP17b}] \label{thm:subOptimality}
Let $\delta_0$ be given by Theorem \ref{thm:well_posed_feedback}.
There exists $\delta \in (0,\delta_0]$ and a constant $C>0$ such that for all $y_0 \in Y$ with $\| y_0 \|_Y \leq \delta$, the following estimates hold:
\begin{align*}
& \mathcal{J}(\mathbf{U}_p(y_0),y_0) \leq \mathcal{V}(y_0) + C
\|y_0\|_{Y}^{p+1}, \label{eq:subOptimality1} \\
& | \mathcal{V}(y_0)-\mathcal{V}_p(y_0) | \leq C \|y_0\|_{Y}^{p+1}.
\end{align*}
Moreover, for all $y_0 \in Y$ with $\| y_0 \|_Y \leq \delta$, problem \eqref{eq:setup} with initial condition $y_0$ possesses a solution $\bar{u}$ satisfying
\begin{align*}
& \| \bar{u} - \mathbf{U}_p(y_0) \|_{L^2(0,\infty;\R^m)} \leq C \| y_0 \|_Y^{(p+1)/2} \\
& \| S(\bar{u},y_0) - S(\mathbf{u}_p,y_0) \|_{W_\infty} \leq C \| y_0 \|_Y^{(p+1)/2}.
\end{align*}
\end{theorem}

\begin{remark} \label{rem:locality}
The constants $\delta_0$, $\delta$, and $C$ involved in Theorem \ref{thm:well_posed_feedback} and \ref{thm:subOptimality} depend on $p$.
They also depend on the data of the problem. In particular, when $\beta$ converges to 0,
the algebraic Riccati equation \eqref{eq:algebraicRiccati} becomes degenerate and the
operator norm of the right-hand sides of the Lyapunov equations \eqref{eqLyapounov3a} possibly increases,
because of the factor $\frac{1}{2\beta}$. Therefore, one can expect that the radius of convergence
of the Taylor expansion and the constant $\delta_0$ both converge to 0 as $\beta$ converges to 0.
\end{remark}

\section{Fokker-Planck equation}
\label{sec:FokkerPlanck}

We describe in this section a specific optimal control problem of the form \eqref{eq:setup} which we shall investigate numerically in
Section \ref{sec:num_res}.

\subsection{Problem formulation}

Following the setup discussed in \cite{BreKP17a}, we consider the following
controlled Fokker-Planck equation:
 \begin{equation}\label{eq:FP_setup}
  \begin{aligned}
  \frac{\partial \rho}{\partial t} &=  \nu \Delta \rho +
\nabla \cdot (\rho \nabla G)+ \sum_{j=1}^m u_{j} \nabla \cdot (\rho
\nabla \alpha_{j})  && \text{in
} \Omega \times (0,\infty), \\
%%%%
0 &= ( \nu \nabla \rho + \rho \nabla G) \cdot \vec{n} && \text{on }
\Gamma \times
(0,\infty), \\
%%%%
\rho(x,0) &= \rho_0(x) && \text{in } \Gamma,
  \end{aligned}
\end{equation}
where $\Omega \in \mathbb R^d $ denotes a bounded domain with
smooth boundary $\Gamma$. The Fokker-Planck equation models the evolution of the
probability distribution of a very large set of particles. More precisely,
$\rho(\cdot,t)$ is the probability density function of the random variable
$X_t$,
solution to the following stochastic differential equation:
\begin{equation*}
\text{d}X(t)= - \nabla_x V(X(t),t) \text{d}t + \sqrt{2\nu} \text{d}W_t,
\end{equation*}
where $(W_t)_{t \geq 0}$ is a Brownian motion and
where the potential $V$ is controlled by $u$ in the following manner:
\begin{equation*}
V(x,t)= G(x) + \sum_{j=1}^m u_j(t) \alpha_j(x), \quad \forall x \in \Omega,\ \forall t \geq 0.
\end{equation*}
Each particle moves along the negative  direction of the gradient  of the potential $V$ and is subject to random perturbations.
When no control is used (i.e.\@ $u=0$), the potential $V$  equals  the \emph{ground potential}  $G$. The functions $\alpha_1$,...$\alpha_m$ are called \emph{control shape functions}.
The reflecting boundary conditions models the fact that the particles are confined in $\Omega$ and
ensure a preservation of probability, i.e.\@
$\int_{\Omega} \rho(x,t) \text{d}x = \int_{\Omega} \rho_0(x) \text{d}x$ for a.e.\@ $t \geq 0.$

The initial probability distribution $\rho_0(x)$ is
normalized so that $\int_{\Omega} \rho_0(x) \mathrm{d}x=1.$
We also assume that the ground potential $G$ and the control shape functions $\alpha_j$ lie in $ W^{1,\infty}(\Omega) \cap
W^{2,\max(2,n)}(\Omega)$ and that $\nabla \alpha_j \cdot \vec{n}=0$ on $\Gamma.$

We introduce now the stationary probability distribution $\rho_\infty$, defined by
\begin{equation*}
\rho_\infty(x)= \frac{e^{-\Phi(x)}}{\int_{\Omega} e^{-\Phi(z)} \text{d} z},
\end{equation*}
where $\Phi(x)= \log(\nu) + \frac{G(x)}{\nu}$. System \eqref{eq:FP_setup} is known to converge to $\rho_\infty$ when $t \to \infty$. This convergence  depends on $\nu$ and the ground potential $G$ and can be extremely slow.
We therefore consider the following optimization problem:
\begin{align}\label{eq:cost_func_FP}
\inf_{u \in L^2(0,\infty;\mathbb R^m)}  \mathcal{J}(u,\rho_0) = \frac{1}{2}
\int_0^\infty \| \rho(t)-\rho_\infty\|
_{L^2(\Omega)}\; \mathrm{d}t + \frac{\beta}{2} \int_0^\infty \| u(t) \|
_{\mathbb R^m}^2 \; \mathrm{d}t,
\end{align}
in order to speed up the convergence to $\rho_\infty$.

\subsection{Abstract formulation and projection}

As is discussed in detail in \cite{BreKP17a} (for the case $m=1$), system
\eqref{eq:FP_setup} can be considered as an abstract bilinear control system of
the form
\begin{equation}\label{eq:abs_pur_bil}
\dot{\rho}(t) = A \rho(t) + \sum_{j=1}^m N_j\rho(t) u_j(t), \quad
\rho(0) = \rho_0,
\end{equation}
where the operators $A$ and $N_j$ are given by
\begin{equation}\label{eq:A_N_op}
\begin{aligned}
  A\colon \mathcal{D}(A)&\subset L^2(\Omega) \to
L^2(\Omega),\\
\mathcal{D}(A)&= \left\{\rho \in H^2(\Omega) \left| (\nu \nabla \rho +
\rho \nabla G) \cdot \vec{n}  =0 \text{ on } \Gamma \right. \right\}, \\
A\rho& = \nu \Delta \rho + \nabla \cdot (\rho \nabla G), \\[1ex]
%%%
N_j\colon H^1(\Omega)& \to L^2(\Omega),\ \ N_j\rho =
\nabla \cdot (\rho \nabla \alpha_j),
\end{aligned}
\end{equation}
and where their
$L^2(\Omega)$-adjoints are given by
\begin{equation}\label{eq:A_N_adj_op}
\begin{aligned}
  A^*\colon \mathcal{D}(A^*)&\subset L^2(\Omega) \to
L^2(\Omega),\\
\mathcal{D}(A^*)&= \left\{\varphi \in H^2(\Omega) \left| (\nu \nabla
\varphi ) \cdot \vec{n}  =0 \text{ on } \Gamma \right. \right\}, \\
A^*\varphi & = \nu \Delta \varphi - \nabla G \cdot \nabla \varphi,
\\[1ex]
%%%
N_j^*\colon H^1(\Omega) &\to
L^2(\Omega),\ \ N_j^*\varphi = -\nabla \varphi \cdot \nabla \alpha_j.
\end{aligned}
\end{equation}
Setting $y= \rho-\rho_\infty$, \eqref{eq:abs_pur_bil} is equivalent to
\begin{equation} \label{eq:abstract_sys}
\dot{y}(t) = A y(t) + \sum_{j=1}^m (N_j y(t) + B_j) u_j(t), \quad
y(0) = \rho_0-\rho_\infty,
\end{equation}
where
\begin{equation*}
B_j\colon \R \rightarrow L^2(\Omega), \quad B_jc= c N_j \rho_\infty.
\end{equation*}

Denoting by $\mathbbm{1}$ the constant function on $\Omega$ equal to 1, one can easily see that:
\begin{equation*}
\mathbbm{1} \in \text{im}(A)^\perp, \quad
\mathbbm{1} \in \text{im}(N_j)^\perp, \quad
\mathbbm{1} \in \text{im}(B_j)^\perp.
\end{equation*}
The mass conservation property follows directly from this observation: $\int_{\Omega} y(\cdot,t) \text{d} x = \int_{\Omega} y(\cdot,0) \text{d} x$, $\forall t \geq 0$. Consider the space
\begin{equation*}
Y_P= \Big\{ v \in L^2(\Omega) : \int_{\Omega} v \text{d} x = 0 \Big\}.
\end{equation*}
The mass conservation property implies that $y(t)$ does not converge to 0 if $y(0)$ does not lie in $Y_P$. Therefore, condition $(A4)$ is not satisfied if \eqref{eq:abstract_sys} is considered as a dynamical system in $L^2(\Omega)$.
Instead, it must be regarded as a dynamical system in $Y_P$. As detailed in
\cite{BreKP17a}, this can be done by first considering the projection $P$ on
$\mathbbm{1}^\perp$ along $\rho_\infty$:
\begin{equation*}
P\colon L^2(\Omega) \rightarrow Y_P, \quad Py= y- \Big( \int_{\Omega} y
\text{d} x \Big) \rho_\infty.
\end{equation*}
Then, we have:
\begin{equation} \label{eq:proj_sys}
\dot{y}(t)= \widehat{A}y(t) + \sum_{j=1}^m (\widehat{N}_j y(t)+\widehat{B}_j)u(t),
\end{equation}
where
\begin{align*}
\widehat{A}= \ & PAI_P \quad \text{with }\mathcal{D}(\widehat{A})= \mathcal{D}(A) \cap Y_P, \\
\widehat{N}_j= \ & PN_j I_P \quad \text{with }\mathcal{D}(\widehat{N}_j)= H^1(\Omega) \cap Y_P, \\
\widehat{B}_j= \ & PB_j,
\end{align*}
and where $I_P$ denotes the injection of $Y_P$ into $L^2(\Omega)$. Assumptions (A1)-(A4) are now satisfied for system \eqref{eq:proj_sys}, as proved in \cite[Section 8]{BreKP17b}.

\section{Algorithmic approach}
\label{sec:alg_app}

Our numerical implementation of the feedback laws is based on the following
approach. We first discretize the Fokker-Planck equation with a
finite-difference scheme, leading to a finite-dimensional bilinear optimal
control problem. Because of the curse of dimensionality, the tensors
$\mathcal{T}_2, \mathcal{T}_3,\dots$ cannot be directly computed for the
discretized problem. A reduction of the discretized model is therefore
necessary. The Lyapunov equations \eqref{eqLyapounov3} can then  be solved
using techniques from \cite{Gra04}.

\subsection{Discretization}

\subsubsection{Discretization of the original state equation}

The spatial discretization is obtained with a finite-difference method.
We use a uniform grid with $n$ points. The
discrete approximations of the individual operators are subsequently denoted
with a subscript $n$. Due to the simpler structure of the boundary conditions
for the operator $A^*$, we employ a finite-difference
scheme for $A^*$ rather than for $A$ itself. Then, the transpose of the resulting matrix serves as
a discrete approximation of $A$. It is denoted by $A_n \in \R^{n \times n}$. For the discretization of the
advective term $\nabla G \cdot \nabla \cdot,$ an upwinding-like scheme which
utilizes backward/forward differences based on the sign of the derivatives
$G_{x_i}$ is used. Since the sign of the controls are not known a
priori, central differences are used for the discretization of $N_j$. The
discretization of $N_j$ is denoted by $N_{j,n} \in \R^{n \times n}$.

\begin{remark}
A finite-difference scheme has been used because it preserves the bilinear structure of the Fokker-Planck equation and therefore allows the computation of a reduced-order model and provides a natural way of discretizing the Riccati equation \eqref{eq:algebraicRiccati} and the Lyapunov equations \eqref{eqLyapounov3}.
Other popular schemes for the discretization of the Fokker-Planck equation, like
the Cooper-Chang algorithm \cite{CooC70} or semi-Lagrangian methods (see the
detailed bibliography of \cite{CarS17}) have nice features (in particular,
positivity preservation), however, they do not maintain the bilinear structure.
\end{remark}

\begin{remark}
The central finite-difference scheme used for the operators $N_j$ worked well for our simulations, even though in principle, it might lead to numerical instabilities.
The design of a scheme for discretizing $N_j$ without prior knowledge of the properties of the control 
is still a challenging issue.
\end{remark}

Since the operator $A$ is known to have a real spectrum with $\rho_\infty$ corresponding
to the smallest eigenvalue (in magnitude) (see \cite[Section 3]{BreKP17a}), a discretization $\rho_{\infty,n} \in \R^n$ can be
efficiently computed, even for large scale problems, by an inverse iteration applied to $A_n$.
Denoting by $h_{x_i}$ the mesh size in the direction $x_i$ and setting $\bar{h}
=\prod_{i=1}^d h_{x_i}$ and $\mathbf{1}=[1,\dots,1]^T  \in \mathbb R^n$,
we normalize the stationary distribution $\rho_{\infty,n}$ so that $\bar{h}\mathbf{1}^T \rho_{\infty,n}=1.$
The initial probability distribution is also normalized:
$\bar{h}\mathbf 1^T \rho_{n,0}  =1.$
Finally, we use $B_{j,n}= N_{j,n} \rho_{\infty,n} \in \R^n$ for the discretization of $B_j$ and set $B_n= [B_{1,n} \hdots B_{m,n}] \in \R^{n \times m}$.

All together the spatially discretized problem reads:
\begin{align}
& \min_{u \in L^2(0,\infty;\R^m)}
J_n(y_{0,n},u) := \frac{1}{2} \int_0^\infty \bar{h} \| y_n(t) \|_{\R^n}^2 \text{d}t
+ \frac{\beta}{2} \int_0^\infty \| u(t) \|_{\R^m}^2 \text{d} t, \label{eq:disc_cost} \\
& \qquad \text{subject to:} \quad
\begin{cases}
\begin{array}{l}
\dot{y}_n(t) = A_ny_n(t) + \Big( \sum_{j=1}^m N_{j,n} y_n(t) u_j(t) \Big) + B_n u(t), \\
y(0)= y_{0,n} := \rho_{0,n}-\rho_{\infty,n}. \label{eq:disc_state_eq}
\end{array}
\end{cases}
\end{align}

\subsubsection{Discretization of the projected state equation}

As explained in Section \ref{sec:FokkerPlanck}, the state equation must be
regarded on a subspace $Y_P$ of $L^2(\Omega)$ to guarantee  stabilizability.
The underlying projection must be numerically implemented to allow an efficient
resolution of the algebraic Riccati equation. We recall the main steps of the
computation of the corresponding discretized and projected operators, details
can be found in \cite{BreKP17a}.
Consider the matrix $R \in \R^{n \times n}$ and its inverse, given by
\begin{equation*}
R= \left( \begin{array}{cccc}
1 & &  & \rho_{\infty,1} \\
 & \ddots & & \vdots \\
 & & 1 & \rho_{\infty,n-1} \\
-1 & \hdots & -1 & \rho_{\infty,n}
\end{array}
\right), \quad
%%%%
R^{-1}= \left( \begin{array}{cccc}
1 & & & 0 \\
 & \ddots & & \vdots \\
 & & 1 & 0 \\
1 & \hdots & 1 & 1
\end{array} \right)
-
\left( \begin{array}{ccc}
\rho_{\infty,1} & \hdots & \rho_{\infty,1} \\
\vdots & \vdots & \vdots \\
\rho_{\infty,n-1} & \hdots & \rho_{\infty,n -1} \\
0 & \hdots & 0
\end{array} \right).
\end{equation*}
The $n-1$ first columns of $R$ build a basis of the orthogonal set to the vector $\mathbf{1}$.
We consider the state space
transformation $\begin{bmatrix} \tilde{y}_n(t) \\ z_n(t) \end{bmatrix} =
R^{-1}y_n(t)$, where $\tilde{y}_n(t) \in \R^{n-1}$.
After the state space transformation, we obtain the system
\begin{align} \label{eq:disc_proj_state_eq}
& \begin{bmatrix} \dot{\tilde{y}}_n(t) \\ \dot{z}_n(t) \end{bmatrix} = \big( R^{-1} A_n R \big)
\begin{bmatrix} \tilde{y}_n(t) \\ z_n(t) \end{bmatrix} + \Big( \sum_{j=1}^m \big( R^{-1} N_{j,n}
R \big) \begin{bmatrix}
\tilde{y}_n(t) \\ z_n(t) \end{bmatrix} u_j(t) \Big) +  \big( R^{-1} B_n \big) u(t), \\
& \begin{bmatrix} \tilde{y}_n(0) \\ z_n(0) \end{bmatrix} = R^{-1} y_{0,n}, \notag
\end{align}
where
\begin{align*}
R^{-1}A_n R =\begin{bmatrix} \widetilde{A}_n & 0 \\ 0 & 0 \end{bmatrix}, \quad
R^{-1} N_{j,n} R = \begin{bmatrix} \widetilde{N}_{j,n} & * \\ 0 & 0 \end{bmatrix},
\quad
R^{-1}B_n = \begin{bmatrix} \widetilde{B}_n \\ 0 \end{bmatrix}.
\end{align*}
Because of the normalization of $\rho_{0,n}$ and $\rho_{\infty,n}$, we have
$z_n(0)=0.$ Moreover, the second block row in \eqref{eq:disc_proj_state_eq} is null, therefore $z(t)= 0$ and
\begin{align*}
y(t) = R\begin{bmatrix} \tilde{y}(t) \\ z(t) \end{bmatrix} = R \begin{bmatrix}
\tilde{y}(t) \\ 0 \end{bmatrix} = R \underbrace{\begin{bmatrix} I_{n-1} \\ 0
\end{bmatrix}}_{=:Q \in \mathbb R^{n\times n-1}} \tilde{y}(t).
\end{align*}
Finally, we obtain the following equivalent formulation of problem \eqref{eq:disc_cost}-\eqref{eq:disc_state_eq}:
\begin{align}
& \min_{u \in L^2(0,\infty;\R^m)}
J_n(y_{0,n},u) := \frac{1}{2} \int_0^\infty \| \widetilde{C}_n \tilde{y}_n(t)\|_{\mathbb
R^n} ^2 \text{d}t
+ \frac{\beta}{2} \int_0^\infty \| u(t) \|_{\R^m}^2 \text{d} t, \label{eq:disc_cost_fct} \\
& \qquad \text{subject to:} \quad
\begin{cases} \begin{array}{l}
\dot{\tilde{y}}_n(t) = \widetilde{A}_n \tilde{y}_n(t) + \Big( \sum_{j=1}^m \widetilde{N}_{j,n} \tilde{y}_n(t) u_j(t) \Big) + \widetilde{B}_n u(t), \\
\tilde{y}_n(0)= Q^T R^{-1} y_{0,n},
\end{array}
\end{cases} \label{eq:disc_proj_system}
\end{align}
where $\widetilde{C}_n = \sqrt{\bar{h}} RQ$.

\subsection{Computation of the feedback tensors}

In theory, the polynomial feedback laws associated with the discretized problem \eqref{eq:disc_cost_fct}-\eqref{eq:disc_proj_system} can be obtained by solving the algebraic Riccati equation and the generalized Lyapunov equations associated with the discretized operators $\widetilde{A}$, $\widetilde{N}_1$,...,$\widetilde{N}_m$, $\widetilde{B}$ (for simplicity, we omit the subscript $n$ in this subsection).
However, the generalized Lyapunov equation of order $k$, corresponding to  the discretized system, is equivalent to a linear system of size $(n-1)^k$.  As a remedy, we propose to replace system \eqref{eq:disc_proj_system} by a reduced-order model.
We describe below our approach for reducing the discretized state equation and explain how to solve the corresponding reduced Lyapunov equations.

\subsubsection{Model reduction}
\label{subsec:feedback_tensors}

We construct a reduced-order model for \eqref{eq:disc_proj_system} of the form
\begin{align}\label{eq:rom}
 \dot{y}_r(t)&= A_r y_r(t) +\sum_{j=1}^m N_{j,r} y_r(t) u_j(t) + B_r u(t), \
y_r(0)= y_{0,r},
\end{align}
where the matrices $A_r, N_{j,r} \in \mathbb R^{r\times r}, B_r \in \mathbb R^{r \times m},
r\ll n-1$ are computed in such a way that for some matrix $C_r \in \R^{n \times r}$, $C_r y_r(t) \approx \widetilde{C}
\tilde{y}_n(t)$, for a range of controls $u_j(t).$
Our construction is based on a known generalization of the
method of balanced truncation for bilinear systems, see e.g.\@ \cite{BenD11}.
It has already been used in the context of the Fokker-Planck equation
in \cite{BenBHS17}. Let us briefly summarize it. As in the case
of linear systems, a reduced-order model is obtained as a \textit{truncation}
of a system that is \textit{balanced} with respect to certain Gramians. In the
bilinear case (\cite{BenD11}), reachability and observability of a bilinear
system can be associated with the definiteness of the Gramians $X$ and $Y$ given
as the solution of the generalized Lyapunov equations
\begin{align*}
  \widetilde{A}X + X\widetilde{A}^T + \sum_{j=1}^m \widetilde{N}_j X \widetilde{N}_j^T +
\widetilde{B}\widetilde{B}^T &=0,
\\
  \widetilde{A}^TY + Y\widetilde{A} + \sum_{j=1}^m \widetilde{N}_j^TY\widetilde{N}_j +
\widetilde{C}^T\widetilde{C} &=0.
\end{align*}
Since an explicit computation of $X$ and $Y$ based on vectorization requires
$\mathcal{O}(n^6)$ operations, we use a fixed point iteration, as discussed in
\cite{Dam08}. More precisely, we compute
\begin{align*}
\widetilde{A}X_1 + X_1 \widetilde{A}^T + \widetilde{B}\widetilde{B}^T &= 0, \quad
  \widetilde{A}X_i + X_i \widetilde{A}^T + \sum_{j=1}^m \widetilde{N}_jX_{i-1}
\widetilde{N}_j^T + \widetilde{B}
\widetilde{B}^T =  0
\end{align*}
and stop when the relative residual
\begin{align*}
  \frac{ \left\|\widetilde{A}X_i + X_i \widetilde{A}^T + \sum_{j=1}^m\widetilde{N}_jX_i
\widetilde{N}_j^T + \widetilde{B} \widetilde{B}^T \right\| _F }{\| \widetilde{B} \widetilde{B}^T \| _F }
\end{align*}
falls below a prescribed tolerance $\varepsilon.$ The same procedure is applied
for computing $Y.$ Once (approximations of) the Gramians $X$ and $Y$ have been
computed, the steps for balancing and truncation are the same as in the linear case.
Based on the product of the Cholesky factors $S$ and $R,$ respectively, of the
Gramians $X=S^TS$ and $Y=R^TR$ a singular value decomposition $U\Sigma V=SR^T$
is computed. Using the best rank-$r$ decomposition of $SR^T$ then yields the
final reduced-order model via a Petrov-Galerkin projection
\begin{align*}
 A_r = W_r^T \widetilde{A} V_r, \ \ N_{j,r} = W_r^T \widetilde{N}_j V_r, \ \ B_r =
W_r^T \widetilde{B}, \ \ C_r = \widetilde{C} V_r,
\end{align*}
where $V_r= S^T U_{(:,1:r)}\Sigma_{(1:r,1:r)}^{-\frac{1}{2}}$ and $W_r = R^T
V_{(:,1:r)} \Sigma_{(1:r,1:r)}^{-\frac{1}{2}}.$ The initial condition is obtained as follows:
\begin{equation*}
y_{0,r} =W_r^T \tilde{y}_n(0).
\end{equation*}
Some comments concerning the
reduced-order modeling approach are in order.

\begin{remark}\label{rem:problems_mor}
In contrast to the linear case,
the generalized method of balanced truncation does not exhibit an a priori
error bound. In the next section, we therefore provide several comparisons
between the original and the reduced model and the corresponding feedback laws.
For applicability of MOR techniques, one typically assumes that the
number of inputs and outputs is small. This is clearly not the case for
$\widetilde{C} \in \mathbb R^{n\times n-1}.$ On the other hand, in case at least
the input space is finite-dimensional, analytic control systems are still known
to have rapidly decaying singular values (\cite{Opm10}).  System theoretic model
reduction techniques typically assume that
the initial value is zero, i.e., $\rho_{0,n}=\rho_{\infty}.$ Obviously, this
leads to a trivial stabilization problem for \eqref{eq:disc_proj_system}. For
nonzero initial values, the initialization of the reduced-oder model is not
obvious and might potentially yield a significantly different transient
response.
While the projected initial condition $y_{r,0}= W_r^T\tilde{y}_n(0)$
 might still lead to deviations
between original and reduced-order model, we expect this effect to be
comparably small since for the theoretical results of the feedback law, the
initial value is assumed to be close to the origin.
\end{remark}

\subsubsection{Lyapunov equations}

It is now possible to compute at a higher degree the feedback laws associated
with the following reduced problem:
\begin{equation} \label{eq:disc_red_prob}
\begin{aligned}
& \min_{u \in L^2(0,\infty;\R^m)}
J_r(y_{0,r},u) := \frac{1}{2} \int_0^\infty \| C_r y_r(t) \|_{\R^n}^2 \text{d}t
+ \frac{\beta}{2} \int_0^\infty \| u(t) \|_{\R^m}^2 \text{d} t, \\
& \qquad \text{subject to:} \quad
\begin{cases}
\begin{array}{l}
\dot{y}_r(t) = A_r y_r(t) + \Big( \sum_{j=1}^m N_{j,r} y_r(t) u_j(t) \Big) + B_r u(t), \\
y(0)= y_{0,r}.
\end{array}
\end{cases}
\end{aligned}
\end{equation}
Note that the above problem has a slightly different structure from problem
\eqref{eq:setup}, because of the operator $C_r$, however, only the algebraic
Riccati equation has to be modified. It reads:
\begin{equation*}
A_r^T \Pi_r + \Pi_r A_r - \frac{1}{\beta} \Pi_r B_r B_r^T \Pi_r  + C_r^T C_r = 0.
\end{equation*}
We set: $A_{\Pi,r}= A_r - \frac{1}{\beta} B_r B_r^T \Pi_r$.
For solving the generalized Lyapunov equations, we represent any multilinear form $S\colon (\R^r)^k \rightarrow \R$ by an array in $\R^{r\times ... \times r}$. The associated vectorization is denoted by $\text{vec}(S) \in \R^{r^k}$.
The generalized Lyapunov equation of order $k$ corresponding to \eqref{eq:disc_red_prob} can be formulated as a tensor-structured linear system:
\begin{align}\label{eq:tensor_red}
  \underbrace{  \sum_{i=1}^k (I \otimes \cdots \otimes I \otimes A_{\Pi,r}^T
\otimes I \otimes \cdots \otimes I )}_{=:\mathbf{A}_{k,r}} \mathrm{vec}(T_{k,r}) = \frac{1}{2\beta}
\sum_{j=1}^m\mathrm{vec}(R_{j,k,r}),
\end{align}
where $\otimes$ is the Kronecker product and where $R_{k,j,r}$ is computed with \eqref{eqLyapounov3b}-\eqref{eqLyapounov3c}.

\begin{remark}
Because of the symmetrization operations involved in \eqref{eqLyapounov3c}, the term $R_{k,j,r}$ must be computed as a sum of $(k-1) + \sum_{i=2}^{k-2} \binom{k}{i} \approx 2^k$ terms.
\end{remark}

Note that an explicit computation of the inverse of $\mathbf{A}_{k,r}$
requires $\mathcal{O}(r^{3k})$ operations which would be infeasible even for
moderate reduced dimensions $r.$ However, the specific tensor structure allows
us to approximate the solution to \eqref{eq:tensor_red} by a
quadrature formula. The method is described and analyzed in
\cite{Gra04}. The main idea consists in
combining an explicit integral representation of the inverse $\mathbf{A}_{k,r}^{-1}$
with a separability property of the matrix exponential of tensor-structured
matrices. The obtained approximation of $\mathbf{A}_{k,r}^{-1}$ takes the form
\begin{align*}
 \mathbf{A}_{k,r}^{-1} \approx \sum_{i=-l}^l w_i \bigotimes_{j=1}^k e^{t_i
A_{\Pi,r}^T}.
\end{align*}
We refer to \cite{Gra04} for the choice of the weights and points.
For our numerical simulations, we have used $l= 50$, leading to a sufficiently accurate approximation.

\begin{remark}
Let us emphasize that the model reduction step does not entirely resolve the
\textit{curse of dimensionality}, since the cost of computing the feedback
tensor $T_{k,r}$ still grows exponentially with $k.$ The use of low-rank
tensor formats would possibly allow to increase the degree of the
polynomial approximation of the feedback law. However, this would
introduce a further approximation error. Moreover, in our numerical
examples, we obtained sufficiently accurate approximations of the optimal
control and thus we refrain from a more detailed discussion on tensor
calculus.
\end{remark}

Once the feedback tensors $T_{k,r}$ have been computed up to a degree $p$, we arrive at the following reduced closed-loop system:
\begin{equation} \label{eq:red_closed_loop}
\dot{y}_r(t)= A_r y_r(t) + \Big( \sum_{j=1}^m N_{j,r} y_r(t) \big(
\mathbf{u}_p(y_r(t)) \big)_j \Big) + B_r \mathbf{u}_p(y_r(t)),
\end{equation}
where the reduced feedback law $\mathbf{u}_p$ is given by:
\begin{equation*}
\mathbf{u}_p(y_r)=
-\frac{1}{\beta}\sum_{k=2}^p
\frac{1}{(k-1)!} T_{r,k}(N_r y_r + B_r,y_r,\dots,y_r).
\end{equation*}
If $\| y_{0,r} \|_{\R^r}$ is sufficiently small, then \eqref{eq:red_closed_loop} is well-posed and the feedback law generates a control $u_p \in L^2(0,\infty;\R^m)$, given by
\begin{equation*}
u_p(t)= \mathbf{u}_p(y_r(t)).
\end{equation*}
In the numerical results below, once the control $u_p$ has been computed (by solving the reduced closed-loop system \eqref{eq:red_closed_loop}), its efficiency is tested with the discretized system \eqref{eq:disc_state_eq}, that is to say, by solving:
\begin{equation}
\dot{y}_n(t) = A_ny_n(t) + \Big( \sum_{j=1}^m N_{j,n} y_n(t) \big( u_p(t) \big)_j \Big) + B_n u_p(t), \quad
y_n(0)= y_{0,n}.
\end{equation}

\section{Numerical results}
\label{sec:num_res}

\newlength\fheight
\newlength\fwidth
\setlength\fheight{0.3\linewidth}
\setlength\fwidth{0.4\linewidth}

We report on numerical tests  in dimension 1 and 2,
respectively. The main discussion focuses on the one-dimensional example while
the two-dimensional example should illustrate the applicability of the method
for larger dynamical systems.

All simulations were done on an \intel Xeon(R) CPU E31270 @ 3.40 GHz x 8,
16 GB RAM, Ubuntu Linux 14.04, \matlab \;Version 8.0.0.783
(R2012b) 64-bit (glnxa64). The solutions of the ODE systems are obtained
with the routine \texttt{ode15}. For solving the algebraic Riccati equation,
we use the routine \texttt{care}. The matrix exponential involved in the approximation formula for $\mathbf{A}_{k,r}^{-1}$ is implemented with the routine \texttt{expm}.

\subsection{One-dimensional example}

\begin{figure}[htb]
  \begin{subfigure}[c]{0.5\textwidth}
    \includegraphics{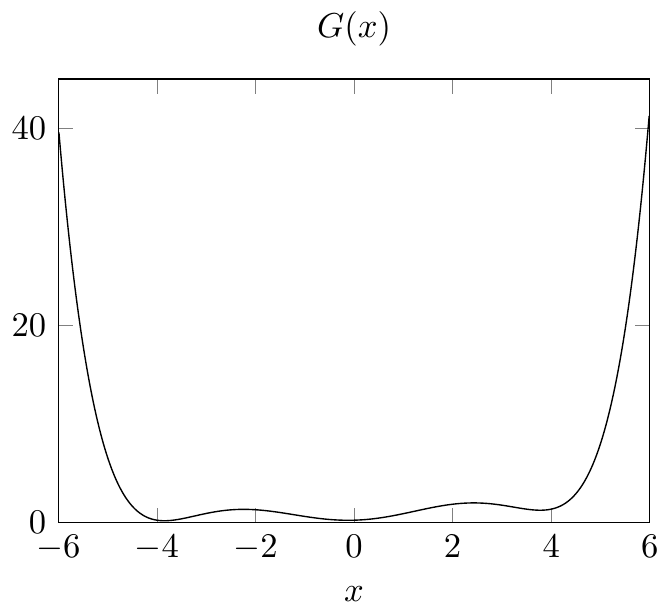}
    \caption{Ground potential.}
    \end{subfigure}\begin{subfigure}[c]{0.5\textwidth}
    \includegraphics{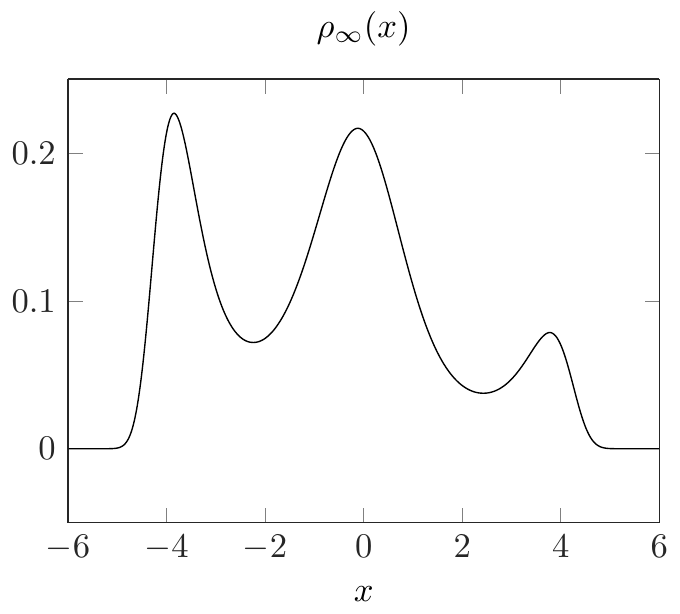}
    \caption{Stationary distribution.}
    \end{subfigure}
  \caption{1D Fokker-Planck equation.}
  \label{fig:1D_FP_setup}
\end{figure}

The first example that we consider is of the form \eqref{eq:FP_setup} with $d=1$, $m=1$,
$\nu=1$ and $\Omega=(-6,6)$.
The ground potential $G$ that we use is represented in Figure \ref{fig:1D_FP_setup}a and the corresponding probability distribution is shown in Figure \ref{fig:1D_FP_setup}b.
The potential $G$ has three local minima reached at $x_1$, $x_2$, and $x_3$ and two local maxima reached at $x_4$ and $x_5$, with
\begin{equation*}
x_1 \approx -3.85 < x_4 \approx -2.24 < x_2 \approx -0.12 < x_5 \approx 2.43 < x_3 \approx 3.78.
\end{equation*}
The minimum is reached at $x_1$. The energy activation $Q$, defined as the highest potential barrier that a particle has to overcome to reach the most stable equilibrium $x_1$, is approximately:
\begin{equation*}
Q= G(y_1)-G(x_2) \approx 1.11.
\end{equation*}
For small values of $\nu$, the rate of convergence of the uncontrolled system is approximately $Ce^{-Q/\nu}$, where $C$ is a constant (see \cite[Section 2]{MatS81}).

The control shape function $\alpha(x) \in \mathbb R$ is such that
\begin{equation*}
\alpha(x) =
\begin{cases}
\begin{array}{cl}
-1/2 & \text{if $-6.0 \leq x \leq -5.9$} \\
x/12 & \text{if $-5.8 \leq x \leq 5.8$} \\
1/2 & \text{if $\phantom{-}5.9 \leq x \leq 6.0$,}
\end{array}
\end{cases}
\end{equation*}
so that $\nabla \alpha \cdot \vec{n}=0$ on $\Gamma$. It is constructed by (twice continuously differentiable) Hermite interpolation on the intervals $(-5.9,-5.8)$ and $(5.8,5.9)$.
The control $u(t)$ is scalar-valued and allows to interact with the
system by tilting one half of the ground potential while raising the other.

%\begin{figure}[htb]
%\begin{center}
%%    \input{controlshape.tikz}
%   \includegraphics{controlshape.pdf}
%   \caption{Control shape function}
%   \label{fig:controlshape}
%  \end{center}
%\end{figure}

Our numerical tests are guided by the following three issues.
First, we show the effect of model reduction on the corresponding feedback laws.
Then, we investigate the convergence of the controls generated by the different
polynomial feedback laws towards the optimal control, as the order $p$ increases.
Finally, we study the influence of the initial condition and the value of $\beta$ on the
efficiency and the convergence of these controls. The last item relates to the local behavior of the method.

\subsubsection{Reduced vs original model}

\begin{figure}[p!]

%% Figure 2
\begin{subfigure}[c]{0.48\textwidth}
\includegraphics{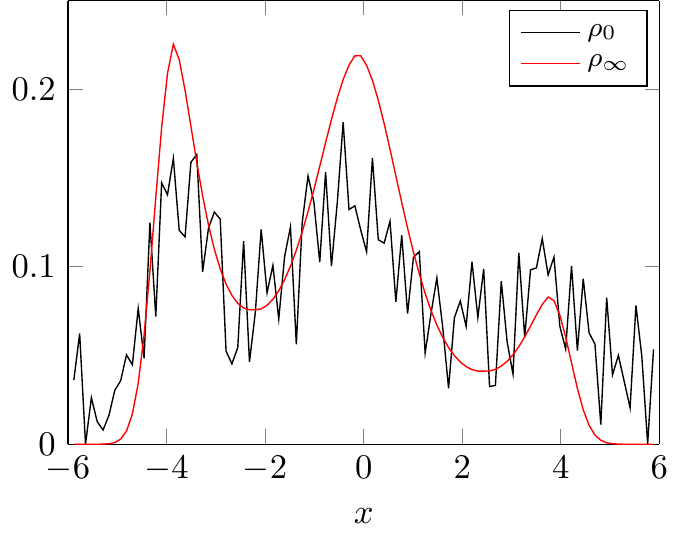}
\caption{Initial and stationary distribution.}
\end{subfigure}
\begin{subfigure}[c]{0.48\textwidth}
\includegraphics{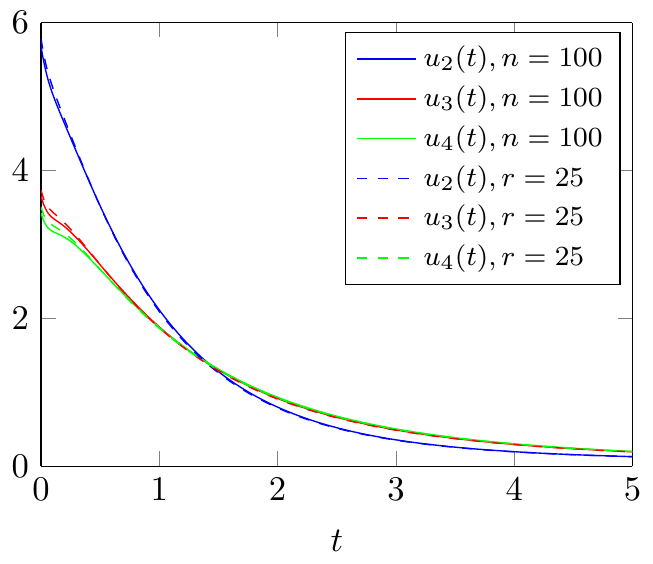}
\caption{Controls.}
\end{subfigure}

\caption{Comparison of the original and reduced models, for $n=100$, $r=25$, and $\beta= 10^{-4}$.}
\label{fig:FOMvsROM_1}

\vspace{10mm}

  \begin{subfigure}[c]{0.5\textwidth}
    \includegraphics{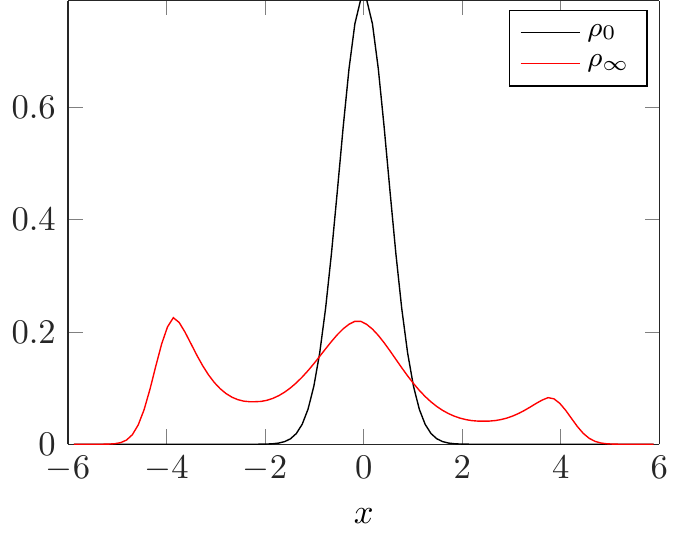}
    \caption{Initial and stationary distribution.}
    \end{subfigure}\begin{subfigure}[c]{0.5\textwidth}
    \includegraphics{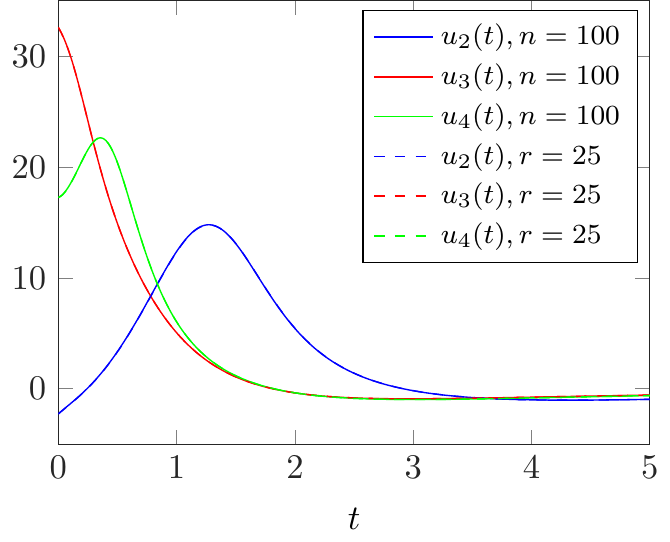}
    \caption{Controls.}
    \end{subfigure}
  \caption{Comparison of the original and reduced models, for $n=100$, $r=25$, and $\beta= 10^{-4}$.}
  \label{fig:FOMvsROM_2}

  \vspace{10mm}

  \begin{subfigure}[c]{0.5\textwidth}
    \includegraphics{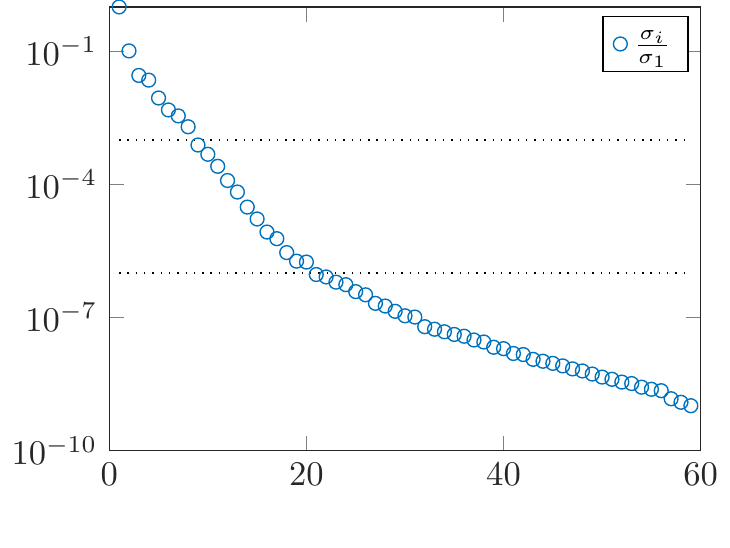}
    \caption{Singular value decay for $n=1000$.}
    \end{subfigure}\begin{subfigure}[c]{0.5\textwidth}
    \includegraphics{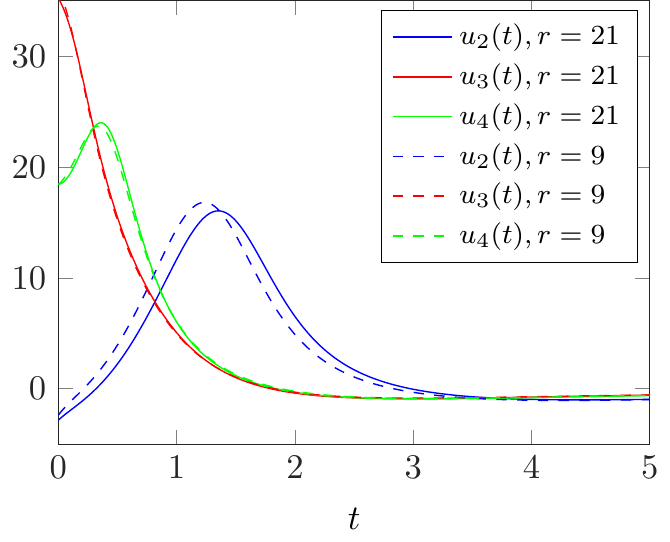}
    \caption{Controls.}
    \end{subfigure}
  \caption{Comparison of the reduced models with $r=21$ and $r= 9$, derived from a finer discretization (with $n= 1000$), with the setup of Figure \ref{fig:FOMvsROM_2}.}
  \label{fig:ROM}

\end{figure}

As described in Subsection \ref{subsec:feedback_tensors}, we rely on a reduced-order model
for the computation of the feedback laws.
We expect the reduced-order model to replicate faithfully the original dynamics.
Due to the absence of a rigorous error bound, we provide the numerical results for some of our test cases.
For this purpose, we first compare the controls obtained with the original model
(for a finite-difference discretization with $n=100$) with the controls
obtained with a reduced model of dimension $r=25$.
The dimension of the
reduced model is determined by neglecting states corresponding to singular
values of the product of the generalized Gramians whose magnitude is smaller
than $10^{-6}.$
Two different initial conditions are considered. The first one is represented in Figure \ref{fig:FOMvsROM_1}a and results from a random perturbation of the stationary distribution.
The second one, represented in Figure \ref{fig:FOMvsROM_2}a, models a set of particles located close to $x=0$, in the second well of the potential $G$.
For both situations, we chose $\beta= 10^{-4}$.
Due to the size of the original model, we are only able to compute the first three
feedback laws when no model reduction is applied.
The controls obtained with the original model and the reduced model are shown in Figures \ref{fig:FOMvsROM_1}b and \ref{fig:FOMvsROM_2}b.
In both cases, the control obtained with the reduced model replicates accurately the ones obtained with the original model.
The only visible deviation occurs at the beginning of the
simulation, for the first initial condition. In view of the nonzero initial condition and Remark
\ref{rem:problems_mor}, this is to be expected.

Since the reduction from $n=100$ to $r=25$ is only moderate, we investigate further
parameters of $n$ and $r,$ respectively. Figure \ref{fig:ROM}a
shows the decay of the singular values of the product of the Gramians for an original
model of dimension $n=1000.$ We include the thresholds for relative magnitudes
smaller than $10^{-3}$ and $10^{-6}$. Let us emphasize that in
contrast to $n=100,$ the relative accuracy of $10^{-6}$ is already obtained for
$r=21$ rather than $r=25.$ For larger values of $n,$ this threshold, however,
remains constant at $r=21.$ Figure \ref{fig:ROM}b shows a comparison
between the controls obtained for reduced-order models of dimension $r=21$
$(\varepsilon=10^{-6})$ and $r=9$ $(\varepsilon=10^{-3}).$ Two comments are in
order: a) comparing Figure \ref{fig:FOMvsROM_2}b with
 Figure \ref{fig:ROM}b, the control laws are visually (almost)
indistinguishable, b) the first singular values remain approximately the same
for discretizations with a larger value of $n$.

\subsubsection{Convergence of higher order feedback laws}

\begin{figure}[htb]
%%% Figure 8
  \begin{subfigure}[c]{0.5\textwidth}
    \includegraphics{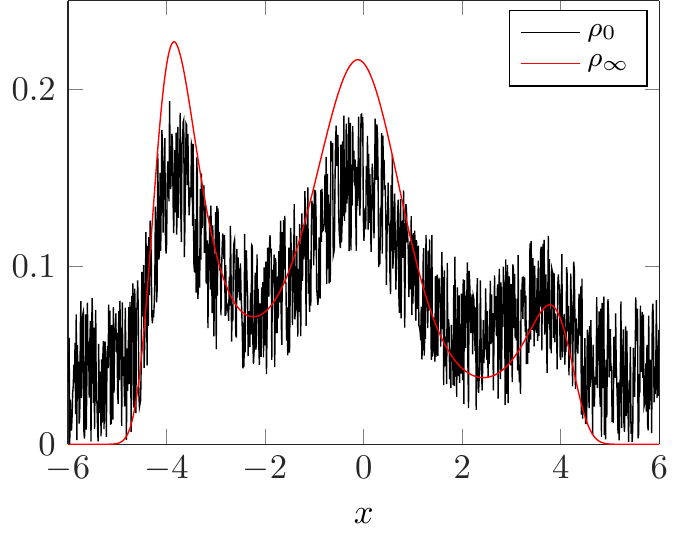}
    \caption{Initial and stationary distribution.}
    \end{subfigure}\begin{subfigure}[c]{0.5\textwidth}
    \includegraphics{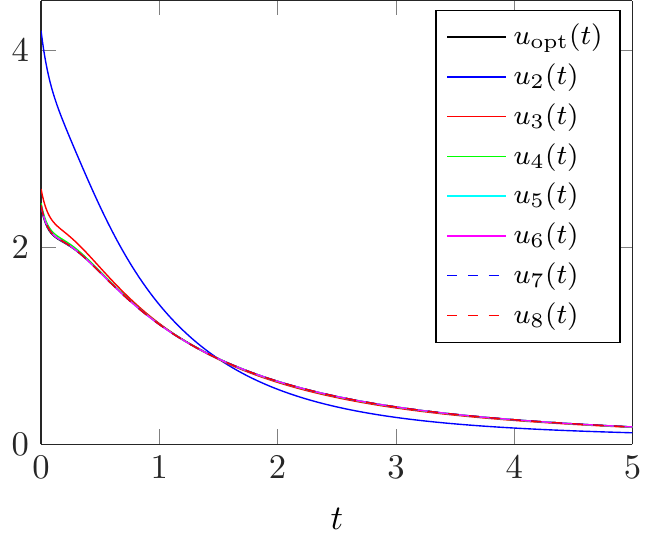}
    \caption{Control laws.}
    \end{subfigure}
  \caption{Convergence of the control laws for $\beta=10^{-4}, n=1000$ and
$r=9.$}
  \label{fig:conv_1}

\vspace{5mm}

%%% Figure 7
  \begin{subfigure}[c]{0.5\textwidth}
    \includegraphics{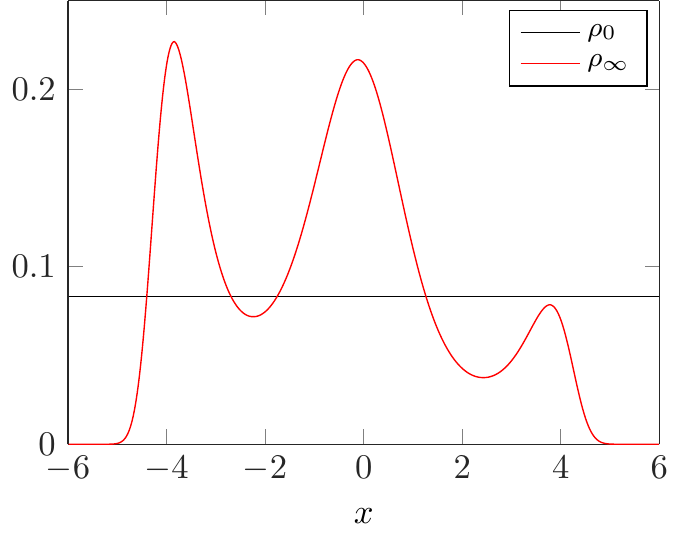}
    \caption{Initial and stationary distribution.}
    \end{subfigure}\begin{subfigure}[c]{0.5\textwidth}
    \includegraphics{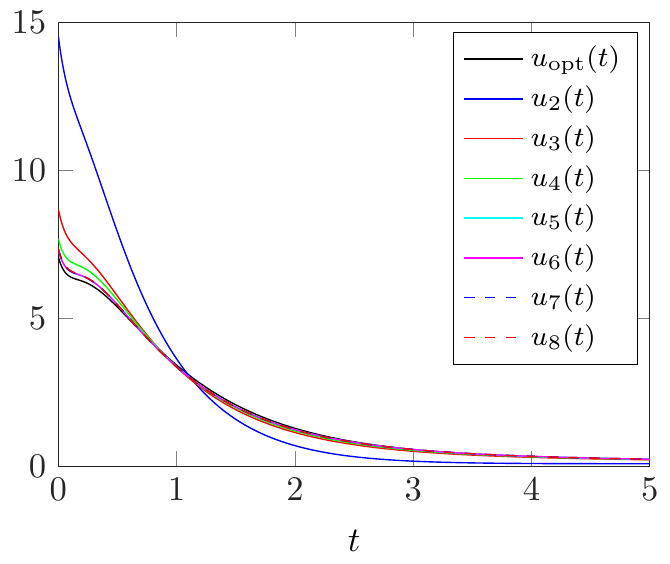}
    \caption{Control laws.}
    \end{subfigure}
  \caption{Convergence of the control laws for $\beta=10^{-4}, n=1000$ and
$r=9.$}
  \label{fig:conv_2}
\end{figure}

We investigate in this subsection the behavior of the controls $u_p$
derived from the feedback laws $\mathbf{u}_p$ for large values of $p$. More precisely, we investigate
the convergence of $u_p$ towards the solution $u_{\text{opt}}$ of the problem, when $p$ increases.
The method  used for computing $u_{\text{opt}}$ is described below.
Two different initial conditions are tested. The first one, represented in
Figure \ref{fig:conv_1}a, is a random perturbation of the stationary distribution. The second one
is the uniform distribution on $\Omega$. We use $\beta= 10^{-4}$ and choose $n= 1000$ for the
discretization. The original model is reduced to $r=9$ $(\varepsilon=10^{-3})$, so that
the seven first feedback laws can be computed. The obtained controls are shown in Figures \ref{fig:conv_1}b and \ref{fig:conv_2}b, respectively.

As can be observed in the case of the randomly perturbed initial condition, the Riccati-based feedback law differs
significantly from all higher order feedback laws. Let us emphasize that the
bilinear term characterized by the operator $N$ does not influence the
computation of the first feedback tensor $\mathcal{T}_2.$ This potentially
explains the strong deviations between $u_2$ and all other controls. We
further see that the higher order control laws quickly approach the optimal
control law $u_{\mathrm{opt}}.$
A clear deviation between the Riccati-based feedback law $u_2$ and all higher order
controls can also be observed for the case of a uniform initial condition. The convergence, however, appears to be slightly slower than in the first case as is indicated by a deviation from $u_3$ and $u_4$ from the other controls. This might be due to the different initial condition which is further
away from the stationary distribution, i.e., $y_0$ is further away from the
origin.

\paragraph{Solving the open-loop problem}

An approximation of the solution of the problem, denoted by $u_{\text{opt}}$, is obtained by
solving
\begin{equation} \label{eq:reduced_open_loop_pb}
\min_{u \in L^2(0,T)} \tilde{J}_r(u):= \frac{1}{2} \int_0^T \| C_r y_r(t) \|_{\R^n}^2 + \beta \| u(t) \|_{\R^m}^2 \text{d}t,
\end{equation}
where $y_r$ is the solution to the reduced-order model \eqref{eq:rom} and where $T= 20$.
To this purpose, we use a gradient-descent algorithm: $u_{k+1}= u_k - t_k \nabla \tilde{J}_r(u_k)$, where
$t_k$ is computed with Armijo's stepsize-rule:
\begin{align*}
t_k= \max_{j=0,1,2,...} \big\{ C \theta^j \,|\,
\tilde{J}_r(u_k - C \theta^j \nabla \tilde{J}_r(u_k) ) \leq \tilde{J}_r(u_k) - C \sigma \theta^j \| \nabla \tilde{J}_r(u_k) \|_{L^2(0,T)}^2 \big\},
\end{align*}
with $C= 500$, $\theta= 0.7$, and $\sigma= 0.05$.
The stopping criterion $\| \nabla \tilde{J}_r(u_k) \|_{L^2(0,T)} \leq \delta$ is used with $\delta= 3\cdot 10^{-4}$. Note that the control provided by such a numerical method may only be an approximation of a local solution to the problem. Even though the gradient-descent algorithm is rather slow, it has the advantage, in the current framework, of being easy to implement and robust. Since the focus of our study is the computation and the analysis of feedback laws, more sophisticated methods for solving \eqref{eq:reduced_open_loop_pb} have not been considered. Let us mention that variants of $\tilde{J}_r$ incorporating a penalty term on the final state $y_r(T)$ provide extremly similar solutions to the problem, since the chosen value for $T$ is large.

\subsubsection{Influence of the initial value and the control costs}

The efficiency of the polynomial feedback laws is only guaranteed in a neighborhood of the origin (i.e.\@ for $\rho_0$ sufficiently close to $\rho_\infty$, in the context of the Fokker-Planck equation). The size of the neighborhood may decrease for small values of $\beta$, as explained in Remark \ref{rem:locality}.
In this subsection, we investigate the efficiency and the convergence of the controls $u_p$ for three different initial conditions, which are respectively close, rather close, and far from the stationary distribution. Different values of $\beta$ are tested. For the discretization and for the reduction of the model, the values $n=1000$ and $r=21$ ($\varepsilon= 10^{-6}$) are used. The integral \eqref{eq:disc_cost_fct} is reduced to the interval $(0,T)$ with $T= 20$ for the evaluation of the cost function.

\paragraph{Test case 1: uniform initial condition}

Table \ref{table:cont_cost_0} provides results for a uniform initial condition. This initial condition can be regarded as very close to the stationary distribution, since the cost of the uncontrolled system $J(u=0)= 0.045$ is small. We have $\| \rho_0- \rho_\infty \|_{L^2(\Omega)} = 0.24$.
In such a situation, the control $u_2$ generated by the feedback law $\mathbf{u}_2$ is almost optimal, as can be seen on Table \ref{table:cont_cost_0}a. For $p=6$, the $L^2$-distance of $u_p$ to $u_{\text{opt}}$  is approximately 7 times smaller than for $p=2$, for the three considered values of $\beta$.

\begin{figure}[htb]
\begin{subfigure}[c]{\textwidth}
\begin{center}
\begin{tabular}{|R{1.4cm}||R{1.4cm}|R{1.4cm}|R{1.4cm}|R{1.4cm}|R{1.4cm}|R{1.4cm}|}
\hline
$\beta$ & $J(u_2)$ & $J(u_3)$ & $J(u_4)$ & $J(u_5)$ & $J(u_6)$ & $J(u_{\text{opt}})$ \\ \hline
1e$^{-3}$ & 0.038 & 0.038 & 0.038 & 0.038 & 0.038 & 0.038 \\
1e$^{-4}$ & 0.034 & 0.033 & 0.033 & 0.033 & 0.033 & 0.032 \\
1e$^{-5}$ & 0.037 & 0.031 & 0.031 & 0.031 & 0.031 & 0.030 \\ \hline
\end{tabular}
\caption{Cost of the controls $u_p$.}
\end{center}
\end{subfigure}\\[3ex]
\begin{subfigure}[c]{\textwidth}
\begin{center}
\begin{tabular}{|R{1.75cm}||R{1.7cm}|R{1.7cm}|R{1.7cm}|R{1.7cm}|R{1.7cm}|R{1.7cm}|}
\hline
\multirow{2}{*}{$\beta$} & \multicolumn{5}{c|}{$\| u_p-u_{\text{opt}} \|_{L^2(0,T)}$} \\
\cline{2-6}
 & $p=2$ & $p=3$ & $p=4$ & $p=5$ & $p=6$  \\ \hline
1e$^{-3}$ & 0.228 & 0.026 & 0.024 & 0.024 & 0.024 \\
1e$^{-4}$ & 4.26 & 1.19 & 0.82 & 0.61 & 0.61 \\
1e$^{-5}$ & 29.8 & 10.3 & 7.91 & 4.70 & 4.05 \\ \hline
\end{tabular}
\caption{$L^2$-distance between the controls $u_p$ and the optimal control $u_{\text{opt}}$.}
\end{center}
\end{subfigure}
\captionof{table}{Convergence results for the test case 1.}
\label{table:cont_cost_0}
\end{figure}

\paragraph{Test case 2: centered initial distribution}

\setcounter{figure}{\value{figure}-1}

\begin{figure}[p!]

\begin{subfigure}[c]{0.48\textwidth}
\includegraphics{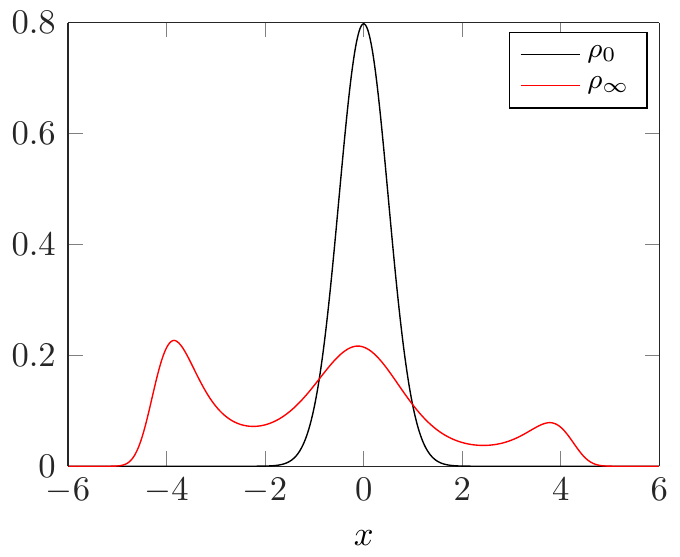}
\caption{Initial and stationary distribution.}
\end{subfigure}
\begin{subfigure}[c]{0.48\textwidth}
\includegraphics{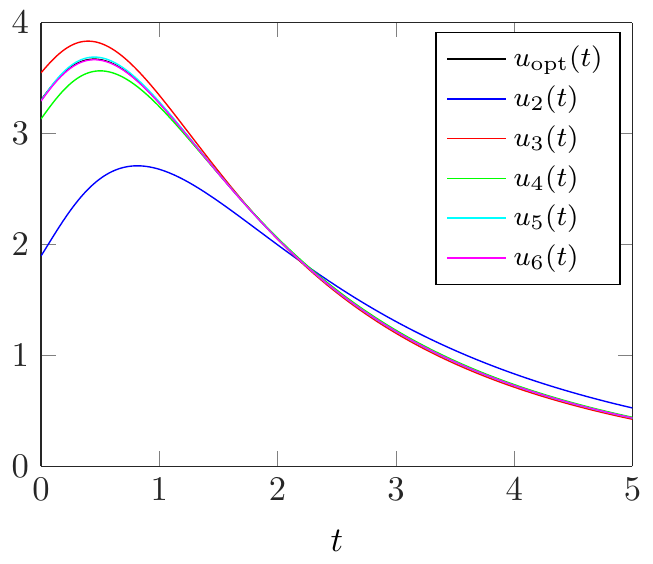}
\caption{Controls for $\beta=10^{-3}$.}
\end{subfigure}\\[5ex]
\begin{subfigure}[c]{0.48\textwidth}
\includegraphics[scale= 0.95]{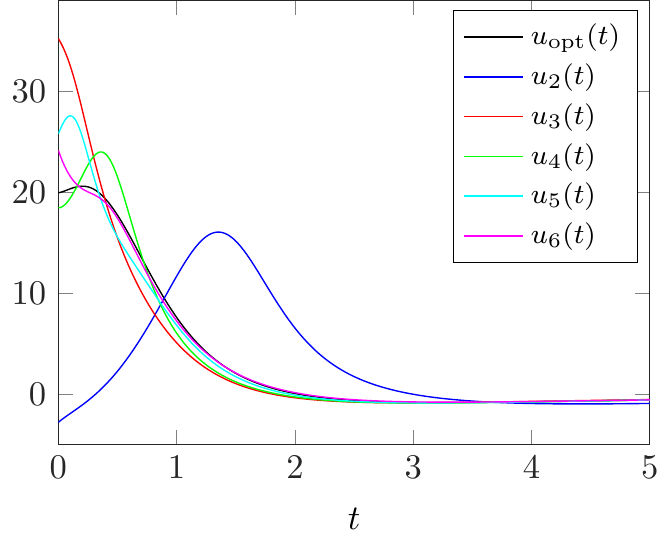}
\caption{Controls for $\beta=10^{-4}$.}
\end{subfigure}
\hspace{-0.6cm}
\begin{subfigure}[c]{0.48\textwidth}
\includegraphics[scale= 0.95]{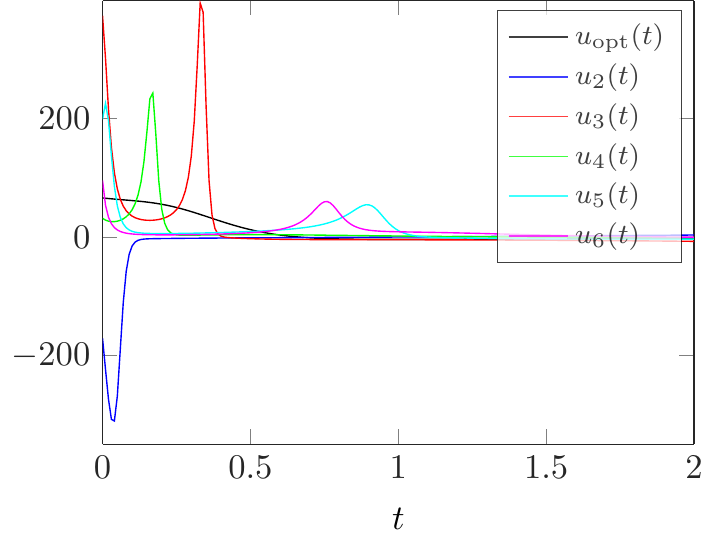}
\caption{Controls for $\beta=10^{-5}$.}
\end{subfigure} \\[2ex]
\caption{Initial condition and controls for the test case 2.}
\label{fig:cont_cost_1}
\vspace{10mm}
\begin{subfigure}[c]{\textwidth}
\begin{center}
\begin{tabular}{|R{1.4cm}||R{1.4cm}|R{1.4cm}|R{1.4cm}|R{1.4cm}|R{1.4cm}|R{1.4cm}|}
\hline
$\beta$ & $J(u_2)$ & $J(u_3)$ & $J(u_4)$ & $J(u_5)$ & $J(u_6)$ & $J(u_{\text{opt}})$ \\ \hline
1e$^{-3}$ & 0.156 & 0.155 & 0.155 & 0.155 & 0.155 & 0.154 \\
5e$^{-4}$ & 0.147 & 0.145 & 0.145 & 0.145 & 0.145 & 0.144 \\
1e$^{-4}$ & 0.138 & 0.122 & 0.120 & 0.120 & 0.120 & 0.119 \\
5e$^{-5}$ & 0.190 & 0.114 & 0.111 & 0.112 & 0.111 & 0.110 \\
1e$^{-5}$ & 0.205 & 0.194 & 0.104 & 0.111 & 0.113 & 0.095 \\ \hline
\end{tabular}
\caption{Cost of the controls $u_p$.}
\end{center}
\end{subfigure} \\[3ex]
\begin{subfigure}[c]{\textwidth}
\begin{center}
\begin{tabular}{|R{1.75cm}||R{1.7cm}|R{1.7cm}|R{1.7cm}|R{1.7cm}|R{1.7cm}|R{1.7cm}|}
\hline
\multirow{2}{*}{$\beta$} & \multicolumn{5}{c|}{$\| u_p-u_{\text{opt}} \|_{L^2(0,T)}$} \\
\cline{2-6}
 & $p=2$ & $p=3$ & $p=4$ & $p=5$ & $p=6$  \\ \hline
1e$^{-3}$ & 1.149 & 0.169 & 0.119 & 0.034 & 0.031 \\
5e$^{-4}$ & 2.583 & 0.737 & 0.171 & 0.336 & 0.219 \\
1e$^{-4}$ & 18.50 & 7.02 & 3.16 & 4.01 & 1.52 \\
5e$^{-5}$ & 46.87 & 13.18 & 8.40 & 8.17 & 2.65 \\
1e$^{-5}$ & 90.5 & 78.0 & 39.0 & 42.6 & 34.3 \\ \hline
\end{tabular}
\caption{$L^2$-distance between the controls $u_p$ and the optimal control $u_{\text{opt}}$.}
\end{center}
\end{subfigure} \\[2ex]
\captionof{table}{Convergence results for the test case 2.}
\label{table:cont_cost_1}
\end{figure}

\setcounter{figure}{\value{figure}-1}

\begin{figure}[p!]
\begin{subfigure}[c]{0.48\textwidth}
\includegraphics{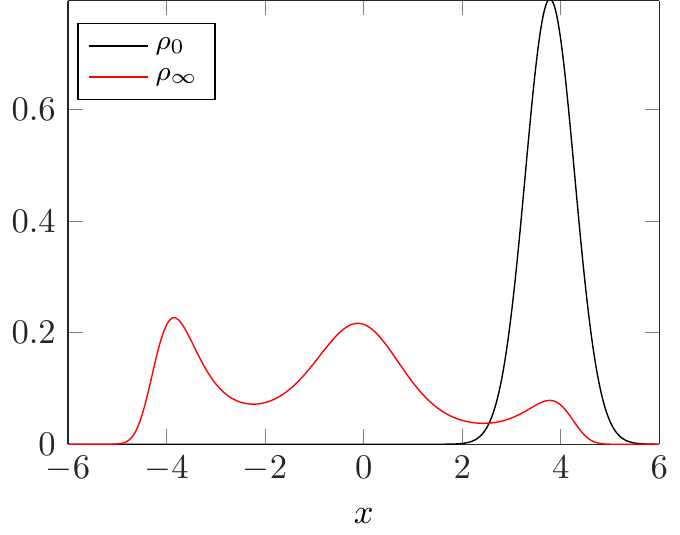}
\caption{Initial and stationary distribution.}
\end{subfigure}
\begin{subfigure}[c]{0.48\textwidth}
\includegraphics{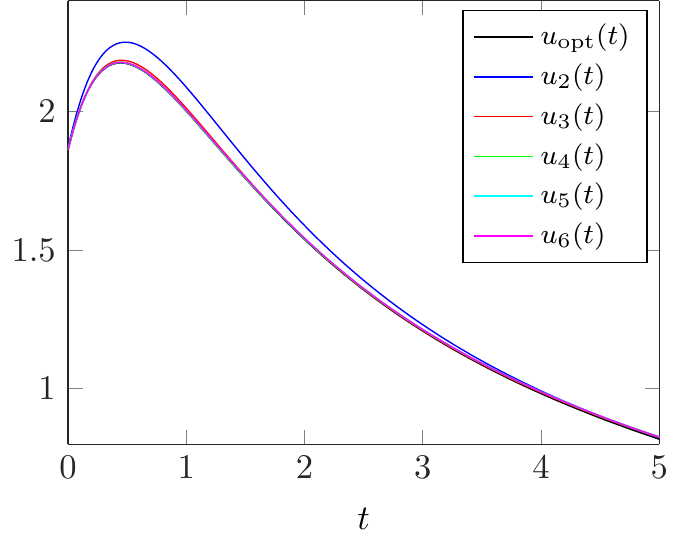}
\caption{Controls for $\beta=10^{-2}$.}
\end{subfigure}\\[5ex]
\begin{subfigure}[c]{0.5\textwidth}
\includegraphics{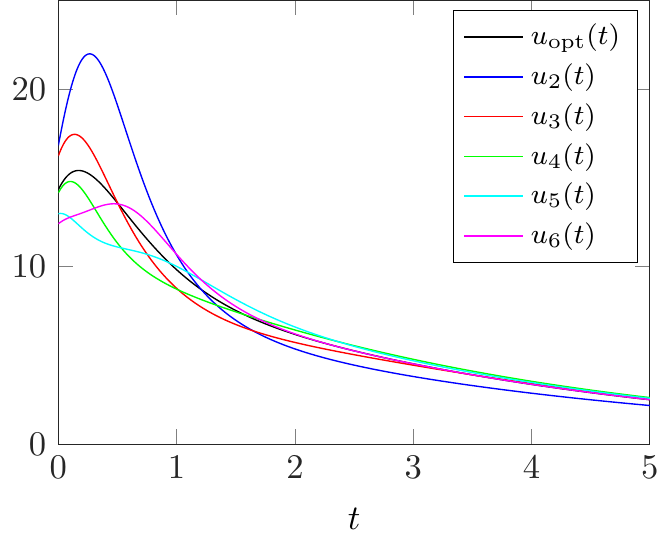}
\caption{Controls for $\beta=10^{-3}$.}
\end{subfigure}\hspace{-0.4cm}
\begin{subfigure}[c]{0.5\textwidth}
\includegraphics{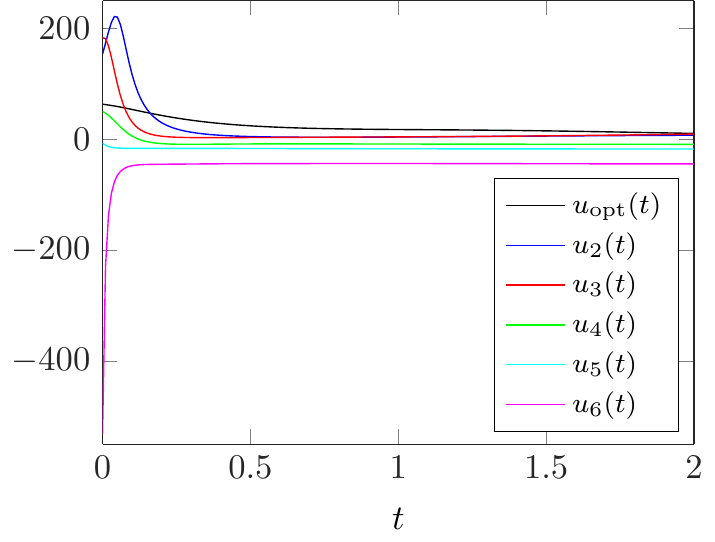}
\caption{Controls for $\beta=10^{-4}$.}
\end{subfigure} \\[2ex]
\caption{Initial condition and controls for the test case 3.}
\label{fig:cont_cost_2}
\vspace{10mm}
\begin{subfigure}[c]{\textwidth}
\begin{center}
\begin{tabular}{|R{1.4cm}||R{1.4cm}|R{1.4cm}|R{1.4cm}|R{1.4cm}|R{1.4cm}|R{1.4cm}|}
\hline
$\beta$ & $J(u_2)$ & $J(u_3)$ & $J(u_4)$ & $J(u_5)$ & $J(u_6)$ & $J(u_{\text{opt}})$ \\ \hline
1e$^{-3}$ & 0.525 & 0.511 & 0.511 & 0.512 & 0.510 & 0.507 \\
5e$^{-4}$ & 0.451 & 0.417 & 0.431 & 0.459 & 0.446 & 0.408 \\
1e$^{-4}$ & 0.381 & 0.368 & 2.689 & $\infty$ & $\infty$ & 0.246 \\
5e$^{-5}$ & 0.381 & 0.432 & $\infty$ & $\infty$ & $\infty$ & 0.206 \\
1e$^{-5}$ & 0.365 & $\infty$ & $\infty$ & $\infty$ & $\infty$ & 0.147 \\ \hline
\end{tabular}
\caption{Cost of the controls $u_p$.}
\end{center}
\end{subfigure}\\[3ex]
\begin{subfigure}[c]{\textwidth}
\begin{center}
\begin{tabular}{|R{1.75cm}||R{1.7cm}|R{1.7cm}|R{1.7cm}|R{1.7cm}|R{1.7cm}|R{1.7cm}|}
\hline
\multirow{2}{*}{$\beta$} & \multicolumn{5}{c|}{$\| u_p-u_{\text{opt}} \|_{L^2(0,T)}$} \\
\cline{2-6}
 & $p=2$ & $p=3$ & $p=4$ & $p=5$ & $p=6$  \\ \hline
1e$^{-3}$ & 4.88 & 1.50 & 1.77 & 2.31 & 1.52 \\
5e$^{-4}$ & 11.26 & 5.03 & 7.11 & 11.89 & 11.99  \\
1e$^{-4}$ & 46.34 & 35.36 & 57.08 & $\infty$ & $\infty$ \\
5e$^{-5}$ & 74.79 & 60.86 & $\infty$ & $\infty$ & $\infty$ \\
1e$^{-5}$ & 172.3 & $\infty$ & $\infty$ & $\infty$ & $\infty$ \\ \hline
\end{tabular}
\caption{$L^2$-distance between the controls $u_p$ and the optimal control $u_{\text{opt}}$.}
\end{center}
\end{subfigure} \\[2ex]
\captionof{table}{Convergence results for the test case 3.}
\label{table:cont_cost_2}
\end{figure}

In this second test case, we consider an initial condition modeling a set of
particles located around the origin. The results are shown on page
\pageref{fig:cont_cost_1} in Figure \ref{fig:cont_cost_1} and Table
\ref{table:cont_cost_1}.
In order to reach the stationary distribution, an important proportion of the particles has to overcome the barrier of the reference potential $G$ located at $y_1 \approx - 2.24$.
The optimal control takes positive values, in order to lower the barrier by tilting the potential on the left
side.
The cost associated with the uncontrolled system is now significantly larger than in the first test case: $J(u=0)= 0.174$. The $L^2$-distance to the equilibrium is larger: $\| \rho_0- \rho_\infty \|_{L^2(\Omega)} = 0.57$.
For all the considered values of $\beta$, a big reduction of $\| u_p-u_{\text{opt}} \|_{L^2(0,T)}$ is observed when the order of the feedback law increases. For $p=6$, the $L^2$-distance is at least 10 times smaller than for $p=2$. Convergence is achieved for values of $\beta$ larger than $5\cdot 10^{-5}$, but is not observed for $\beta= 10^{-5}$, as is well shown in Figure \ref{fig:cont_cost_1}d.
For this intermediate initial condition, the convergence of the controls as well as an important reduction of the costs can be observed, at least for the smallest values of $\beta$.
For values of $\beta$ larger than $5\cdot 10^{-4}$, the controls $u_2$,...$u_6$ are all almost optimal, while for $\beta$ ranging from $10^{-4}$ to $5\cdot 10^{-5}$, a significant difference between $u_2$ and $u_3$ is observed. For $\beta= 10^{-5}$, the costs of $u_4$ are twice smaller as those of $u_2$ and $u_3$.

\paragraph{Test case 3: right-sided initial distribution}

A third test case is presented page \pageref{fig:cont_cost_2} in Figure \ref{fig:cont_cost_2} and Table \ref{table:cont_cost_2}, where the set of
particles is assumed to be located in the third potential well. This initial configuration appears to be more
challenging than the two other configurations, since a large proportion of the set of particles
has now to overcome two barriers of the reference potential, a first one at $y_2 \approx  2.43$ and a second one at $y_1 \approx -2.24$. The cost of the uncontrolled system is $J(u=0)= 0.865$ and the $L^2$-distance of the initial condition to the stationary distribution is $\| \rho_0- \rho_\infty \|_{L^2(\Omega)} = 0.76$
For a comparably high control cost
parameter $\beta=10^{-2},$ the controls rapidly converge. Lowering the
parameter to $\beta=10^{-3},$ convergence is in question (at least cannot
be determined from the numerical results). Finally, for $\beta=10^{-4},$ only
the control laws $u_2$, $u_3$ and $u_4$ actually converge to zero. Higher
order closed-loop system appear to be attracted by a further (nontrivial) steady state.
This is indicated by the symbol $\infty$ in Table \ref{table:cont_cost_2}.
This behavior can be explained by the fact that the closed-loop system is a nonlinear (polynomial) equation
for which different steady states might occur. In the case of the Fokker-Planck equation, the feedback laws $\mathbf{u}_5$ and $\mathbf{u}_6$ introduce a shift of the ground potential such that the particle remains in
the (stable) stationary distribution associated with this shifted potential.
This last test case shows the local nature of the method.

\subsection{A two-dimensional example}

\begin{figure}[p!]
\begin{subfigure}[c]{0.48\textwidth}
\includegraphics[trim= 4cm 9cm 4cm 9cm, clip, scale= 0.5]{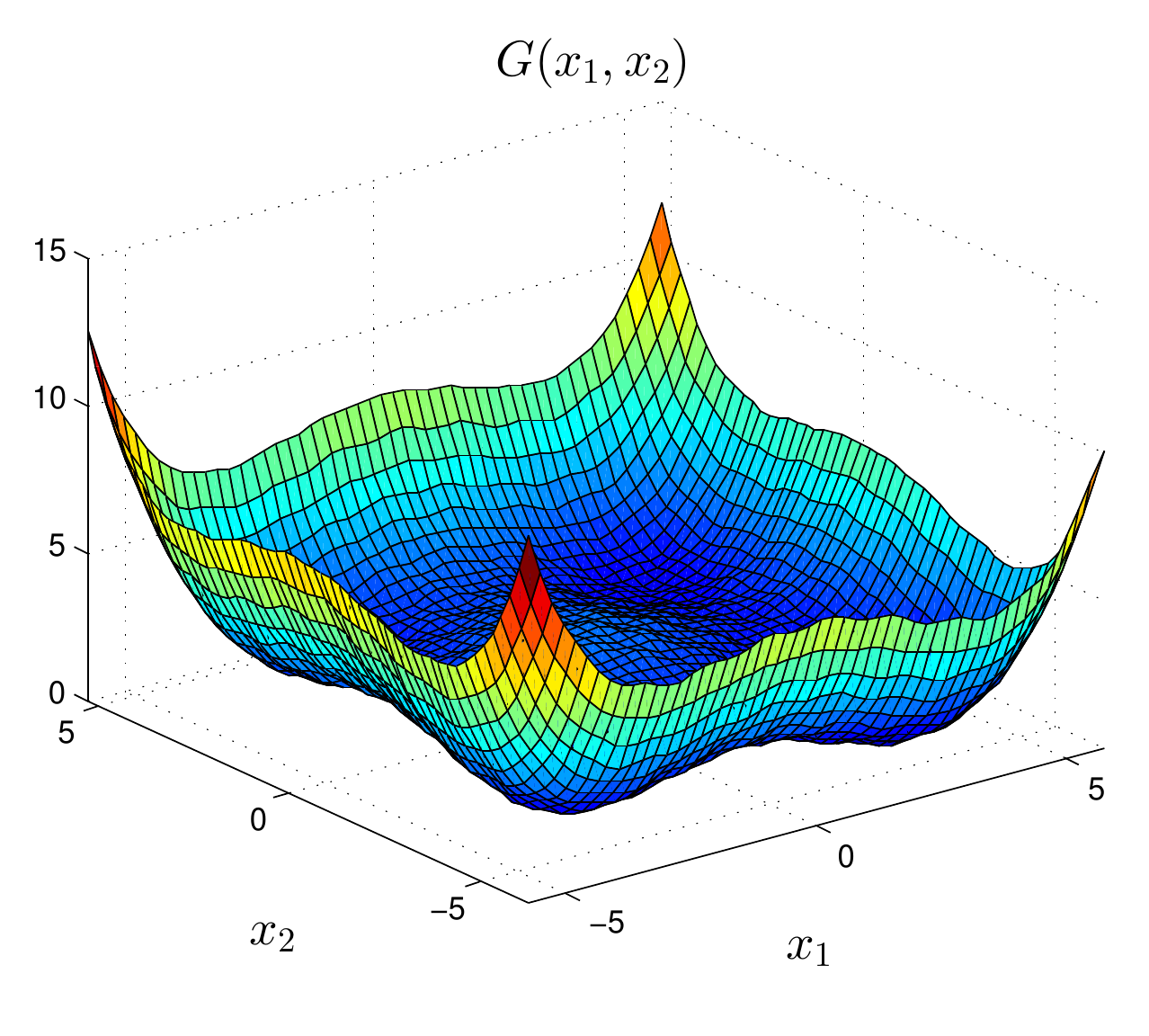}
\caption{Ground potential.}
\end{subfigure}
\begin{subfigure}[c]{0.48\textwidth}
\includegraphics[trim= 4cm 9cm 4cm 9cm, clip, scale= 0.5]{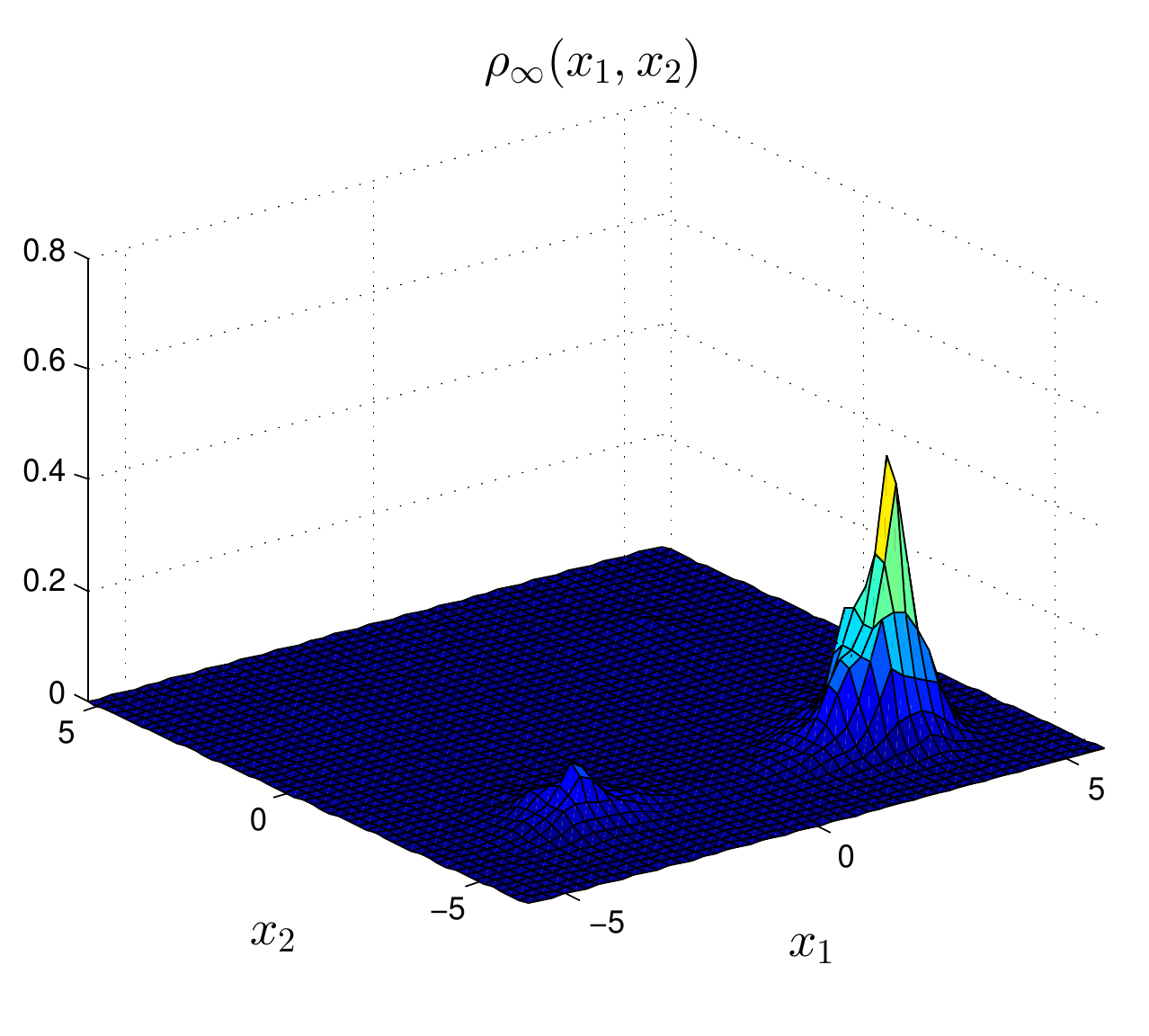}
\caption{Stationary distribution.}
\end{subfigure}
\caption{2D Fokker-Planck equation.}
\label{fig:2d_setup}
\vspace{10mm}
\begin{subfigure}[c]{0.48\textwidth}
\includegraphics[trim= 4cm 9cm 4cm 9cm, clip, scale= 0.5]{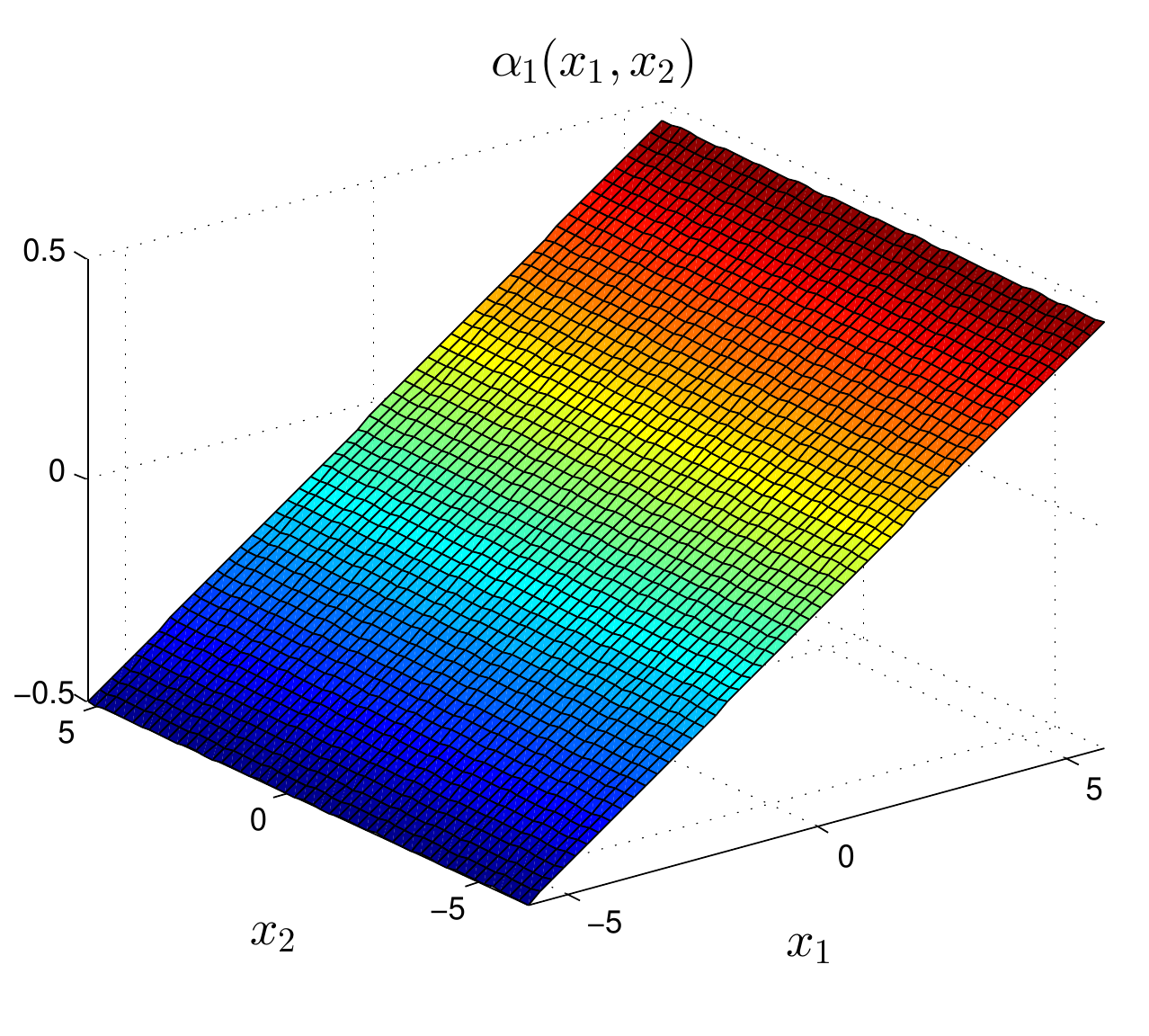}
\end{subfigure}
\begin{subfigure}[c]{0.48\textwidth}
\includegraphics[trim= 4cm 9cm 4cm 9cm, clip, scale= 0.5]{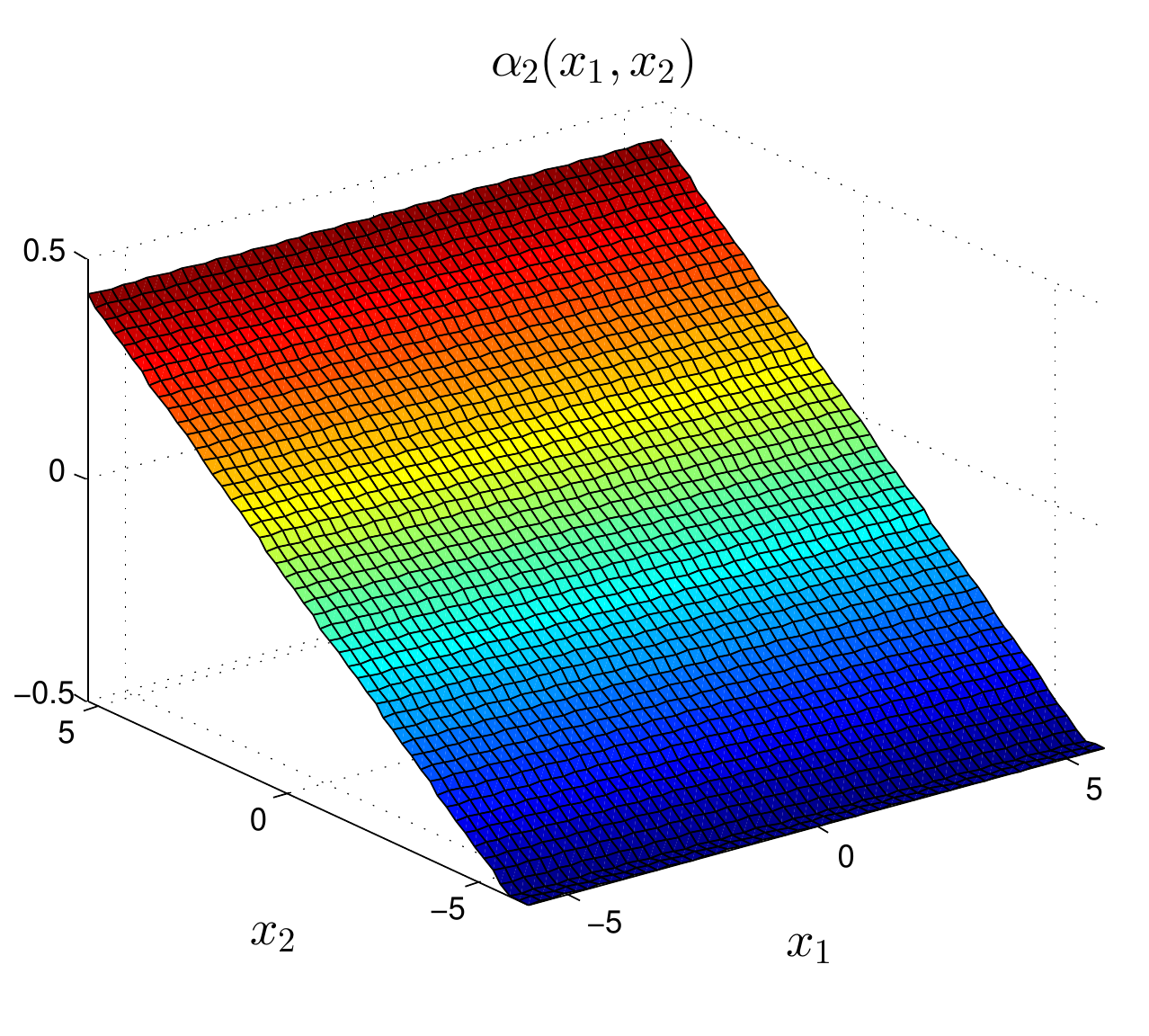}
\end{subfigure}
\caption{Control shape functions $\alpha_1$ and $\alpha_2.$}
\label{fig:cont_shape_2d}
\vspace{8mm}
\begin{center}
\includegraphics{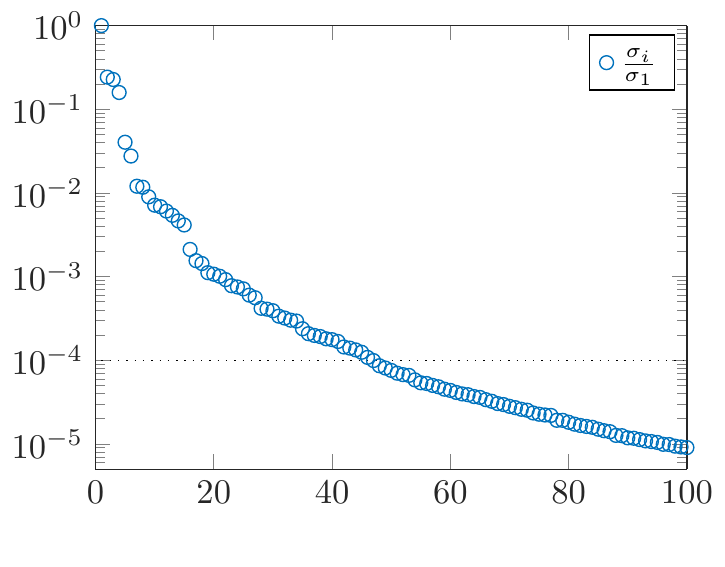}
\caption{Singular value decay for $n=2500$}
\label{fig:value_decay}
\end{center}
\end{figure}

For this second example, we consider a system of the form \eqref{eq:FP_setup} with
\begin{equation*}
d=2, \quad m=2, \quad \nu= 0.25 \quad \text{and} \quad \Omega= (-6,6) \times
(-6.5,5.5).
\end{equation*}
The ground potential $G$ is represented in Figure \ref{fig:2d_setup}a and
the corresponding probability distribution is shown in Figure
\ref{fig:2d_setup}b. The potential $G$ has four local minimizers, located as
follows:
\begin{equation*}
x_A \approx (2.48,-3.77), \quad
x_B \approx (-2.86,-3.75), \quad
x_C \approx (-2.88,2.52), \quad
x_D \approx (2.42,2.51).
\end{equation*}
Two control shape functions are used represented in Figure \ref{fig:cont_shape_2d} and given by:
\begin{align*}
\alpha_1(x_1,x_2)= x_1/12, \quad \forall (x_1,x_2) \in (-5.8,5.8) \times (-6.5,5.5), \\
\alpha_2(x_1,x_2)= x_2/12, \quad \forall (x_1,x_2) \in (-6,6) \times (-6.3,5.3).
\end{align*}
The control shape function $\alpha_1$ is constructed by interpolation on $((-6,-5.8) \cup (5.8,6)) \times (-6,6)$ so that $\nabla \alpha_1 \cdot \vec{n}= 0$ on $\Gamma$, as in the one-dimensional case. The technique is also used for $\alpha_2$.
A negative value of $u_1$ allows to shift the distribution along the first coordinate axis of $\Omega$ and a negative value of $u_2$ allows to shift the distribution along the second coordinate axis.

As for the one-dimensional case, we investigate the influence of the initial condition and the value of $\beta$ on the efficiency of the feedback laws.
We present below the results obtained for two different initial conditions, for
a reduced model of order $r= 47$, obtained from a finite-difference
discretization with $n= 50 \cdot 50$ degrees of freedom with a tolerance of
$\varepsilon= 10^{-4}$. For such a dimension, only the first four feedback
laws can be computed. Figure \ref{fig:value_decay} shows the decay of the
singular values of the product of the Gramians for the unreduced discretized
system. As can be observed, the decay is significantly slower than in the
one-dimensional case. This can be partly explained by the fact that the ground
potential $G$ has a more complicated structure, and that a wider range of
controls are taken into account.
The open-loop control problem is solved with the same parameters:
\begin{equation*}
C=500, \quad \theta= 0.7, \quad \sigma= 0.05, \quad \delta= 3 \cdot 10^{-4}.
\end{equation*}
The final-time used for solving \eqref{eq:reduced_open_loop_pb} and for evaluating \eqref{eq:disc_cost_fct} is set to $T= 200$.

\paragraph{Test case 4: a random perturbation of the initial condition}

\begin{figure}[p!]
\begin{subfigure}[c]{0.48\textwidth}
\begin{center}
\includegraphics[trim= 4cm 9cm 4cm 9cm, clip, scale= 0.5]{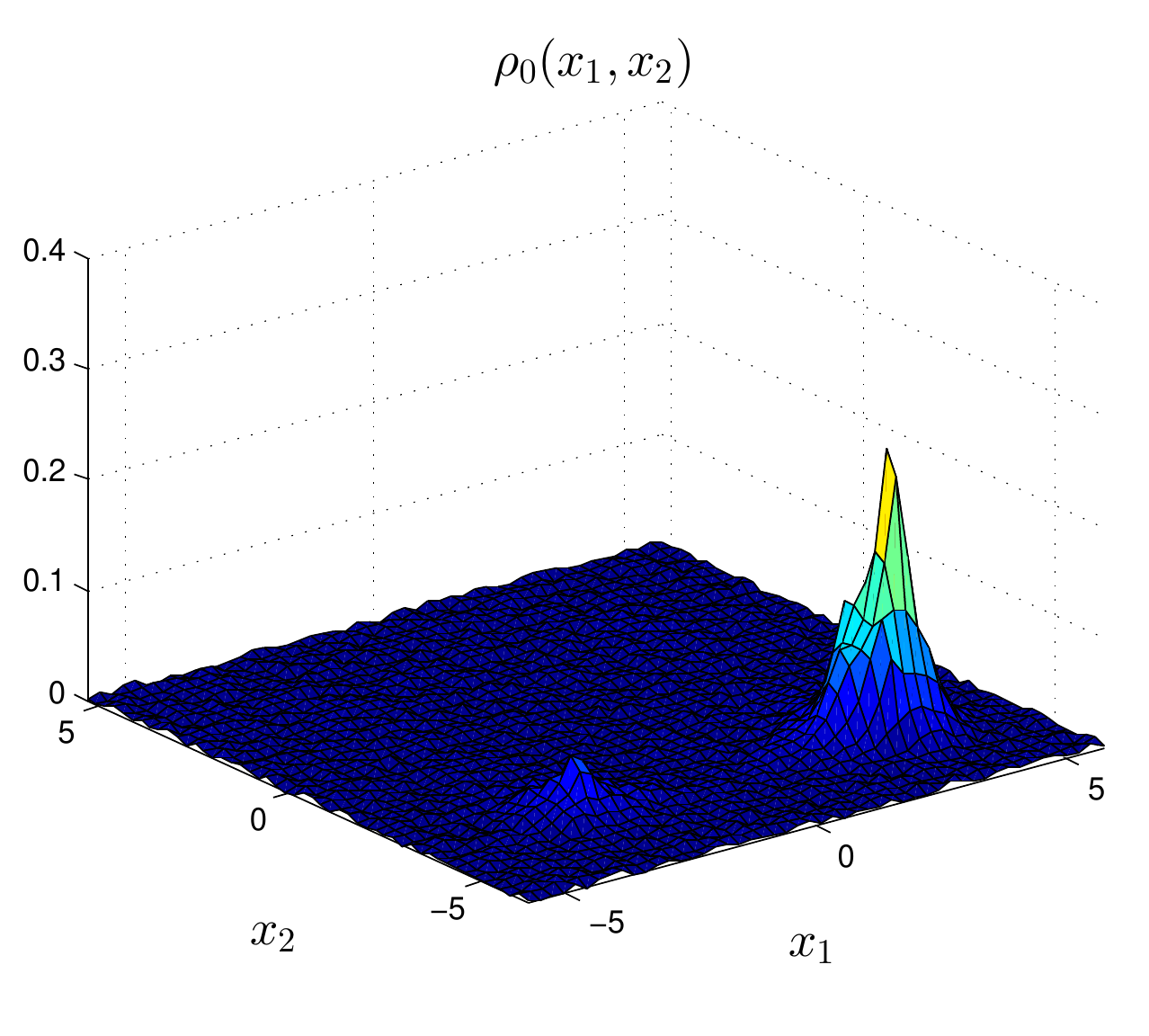}
\caption{Initial distribution.}
\end{center}
\end{subfigure}
\begin{subfigure}[c]{0.48\textwidth}
\begin{center}
\includegraphics{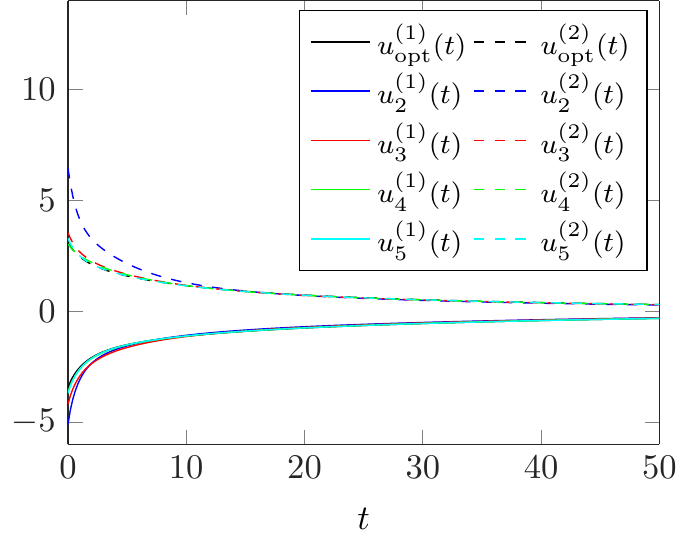}
\caption{Controls for $\beta= 10^{-3}$.}
\end{center}
\end{subfigure}\\[5ex]
\begin{subfigure}[c]{0.48\textwidth}
\begin{center}
\includegraphics{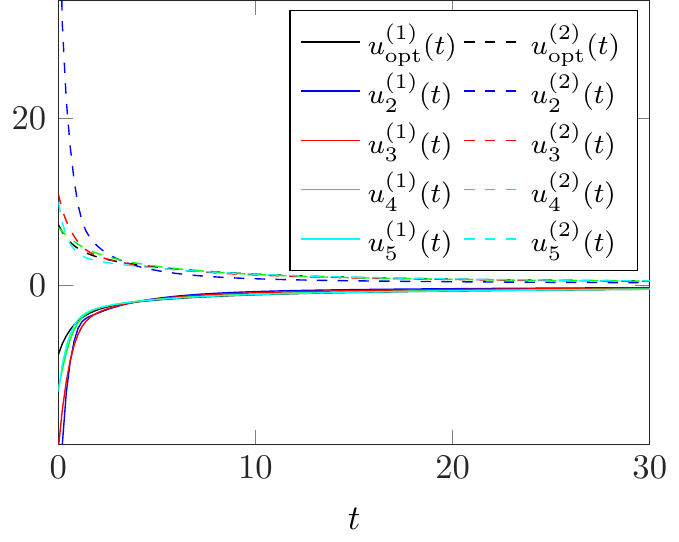}
\caption{Controls for $\beta= 10^{-4}$.}
\end{center}
\end{subfigure}
\begin{subfigure}[c]{0.48\textwidth}
\begin{center}
\includegraphics{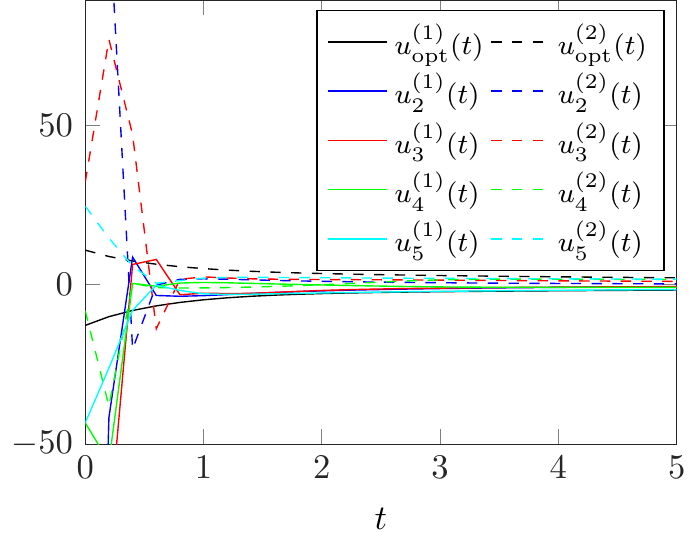}
\caption{Controls for $\beta= 10^{-5}$.}
\end{center}
\end{subfigure} \\[2ex]
\caption{Initial condition and controls for the test case 4.}
\label{fig:2d_random}
\vspace{10mm}
\begin{subfigure}[c]{\textwidth}
\begin{center}
\begin{tabular}{|R{1.35cm}||R{1.35cm}|R{1.35cm}|R{1.35cm}|R{1.35cm}|R{1.35cm}|}
\hline
$\beta$ & $J(u_2)$ & $J(u_3)$ & $J(u_4)$ & $J(u_5)$ & $J(u_{\text{opt}})$ \\ \hline
1e$^{-3}$ & 0.247 & 0.235 & 0.234 & 0.234 & 0.232 \\
5e$^{-4}$ & 0.232 & 0.207 & 0.205 & 0.205 & 0.203 \\
1e$^{-4}$ & 0.252 & 0.180 & 0.174 & 0.174 & 0.171 \\
5e$^{-5}$ & 0.279 & 0.179 & 0.168 & 0.168 & 0.165 \\
1e$^{-5}$ & 0.524 & 0.182 & 20.696 & 0.164 & 0.158 \\ \hline
\end{tabular}
\caption{Cost of the controls $u_p$.}
\end{center}
\end{subfigure}\\[5ex]
\begin{subfigure}[c]{\textwidth}
\begin{center}
\begin{tabular}{|R{1.8cm}||R{1.7cm}|R{1.7cm}|R{1.7cm}|R{1.7cm}|}
\hline
\multirow{2}{*}{$\beta$} & \multicolumn{4}{c|}{$\| u_p-u_{\text{opt}} \|_{L^2(0,T)}$} \\
\cline{2-5}
 & $p=2$ & $p=3$ & $p=4$ & $p=5$ \\ \hline
1e$^{-3}$ & 3.53 & 0.80 & 0.19 & 0.14 \\
5e$^{-4}$ & 6.73 & 1.42 & 0.37 & 0.24 \\
1e$^{-4}$ & 27.40 & 5.78 & 1.83 & 1.24 \\
5e$^{-5}$ & 52.50 & 11.06 & 3.69 & 2.40 \\
1e$^{-5}$ & 257.01 & 63.97 & 84.31 & 10.61 \\ \hline
\end{tabular}
\caption{$L^2$-distance between the controls $u_p$ and the optimal control $u_{\text{opt}}$.}
\end{center}
\end{subfigure} \\[2ex]
\captionof{table}{Convergence results for the test case 4.}
\label{table:2d_random}
\end{figure}

\setcounter{figure}{\value{figure}-1}

\begin{figure}[p!]
\begin{subfigure}[c]{0.48\textwidth}
\begin{center}
\includegraphics[trim= 4cm 9cm 4cm 9cm, clip, scale= 0.5]{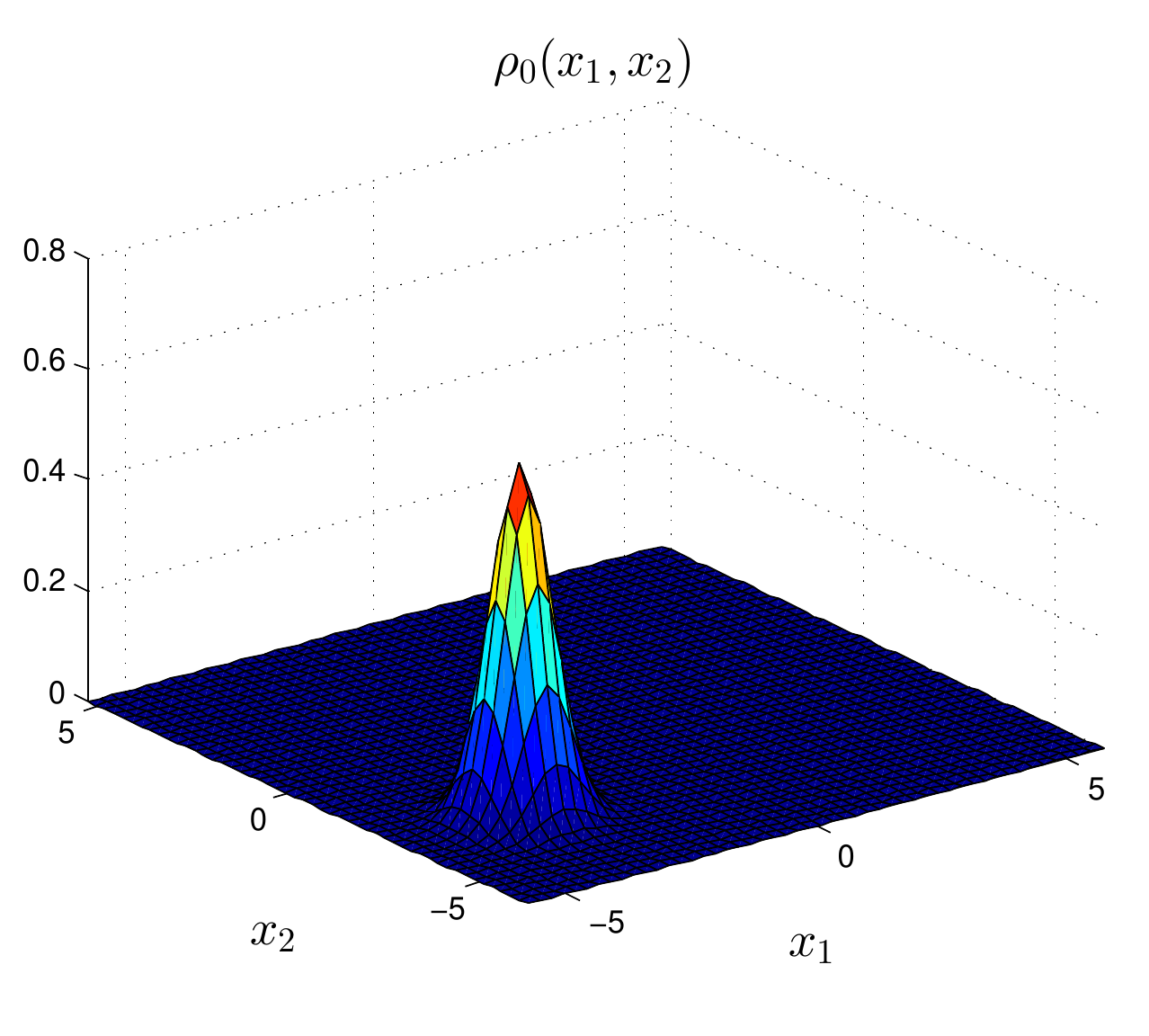}
\caption{Initial distribution.}
\end{center}
\end{subfigure}
\begin{subfigure}[c]{0.48\textwidth}
\begin{center}
\includegraphics[scale=0.9]{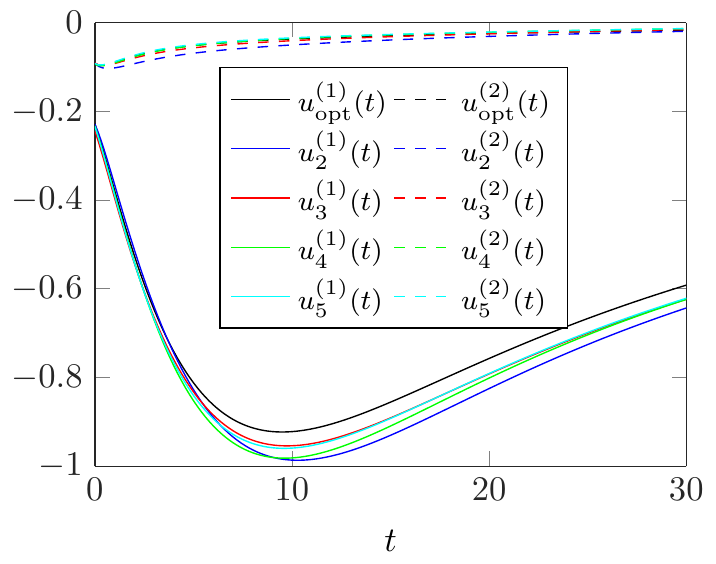}
\caption{Controls for $\beta= 10^{-1}$.}
\end{center}
\end{subfigure}\\[5ex]
\begin{subfigure}[c]{0.48\textwidth}
\begin{center}
\includegraphics[scale=0.9]{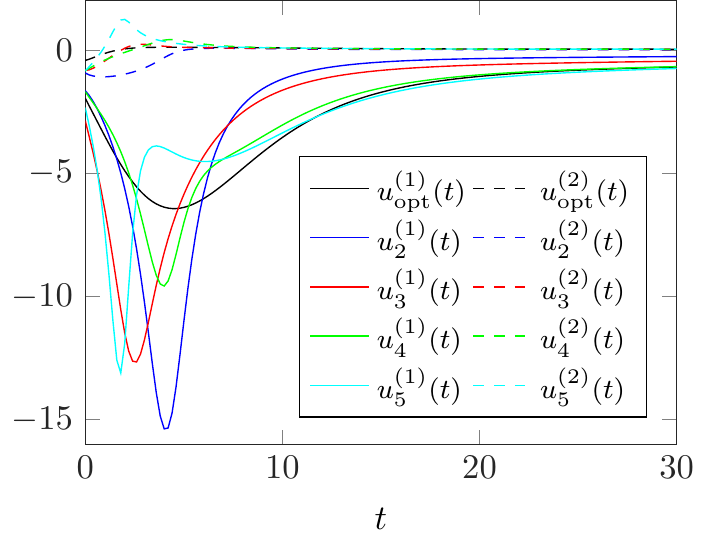}
\caption{Controls for $\beta= 10^{-2}$.}
\end{center}
\end{subfigure}
\begin{subfigure}[c]{0.48\textwidth}
\begin{center}
\includegraphics[scale=0.9]{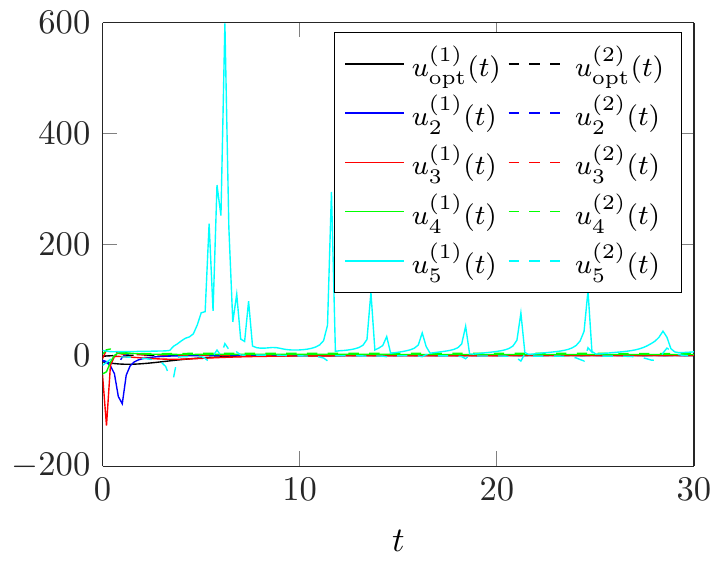}
\caption{Controls for $\beta= 10^{-3}$.}
\end{center}
\end{subfigure} \\[2ex]
\caption{Initial condition and controls for the test case 5.}
\label{fig:2d_gaussian}
\vspace{10mm}
\begin{subfigure}[c]{\textwidth}
\begin{center}
\begin{tabular}{|R{1.35cm}||R{1.35cm}|R{1.35cm}|R{1.35cm}|R{1.35cm}|R{1.35cm}|}
\hline
$\beta$ & $J(u_2)$ & $J(u_3)$ & $J(u_4)$ & $J(u_5)$ & $J(u_{\text{opt}})$ \\ \hline
1e$^{-1}$ & 7.58 & 7.57 & 7.57 & 7.57 & 7.52 \\
5e$^{-2}$ & 6.41 & 6.39 & 6.40 & 6.39 & 6.35 \\
1e$^{-2}$ & 3.70 & 3.34 & 3.09 & 3.32 & 3.00 \\
5e$^{-3}$ & 3.07 & 2.68 & 2.28 & 2.96 & 2.05 \\
1e$^{-3}$ & 2.45 & 2.41 & $\infty$ & $\infty$ & 0.93 \\ \hline
\end{tabular}
\caption{Cost of the controls $u_p$.}
\end{center}
\end{subfigure}\\[5ex]
\begin{subfigure}[c]{\textwidth}
\begin{center}
\begin{tabular}{|R{1.8cm}||R{1.7cm}|R{1.7cm}|R{1.7cm}|R{1.7cm}|}
\hline
\multirow{2}{*}{$\beta$} & \multicolumn{4}{c|}{$\| u_p-u_{\text{opt}} \|_{L^2(0,T)}$} \\
\cline{2-5}
 & $p=2$ & $p=3$ & $p=4$ & $p=5$ \\ \hline
1e$^{-1}$ & 0.70 & 0.61 & 0.62 & 0.60 \\
5e$^{-2}$ & 1.10 & 0.69 & 0.80 & 0.63 \\
1e$^{-2}$ & 13.02 & 11.10 & 4.08 & 9.01 \\
5e$^{-3}$ & 21.59 & 19.80 & 9.66 & 20.06 \\
1e$^{-3}$ & 47.34 & 55.69 & $\infty$ & $\infty$ \\ \hline
\end{tabular}
\caption{$L^2$-distance between the controls $u_p$ and the optimal control $u_{\text{opt}}$.}
\end{center}
\end{subfigure} \\[2ex]
\captionof{table}{Convergence results for the test case 5.}
\label{table:2d_gaussian}
\end{figure}

Figure \ref{fig:2d_random} and Table \ref{table:2d_random} (page \pageref{fig:2d_random}) show the results obtained for an initial condition obtained by randomly perturbing the stationary distribution. The initial condition is therefore close to the stationary distribution. The cost of the uncontrolled system is $J(u=0)= 0.593$ and $\|\rho_0 - \rho_\infty \|_{L^2(\Omega)} = 0.18$. Good convergence results are observed for values of $\beta$ ranging from $10^{-3}$ to $5\cdot 10^{-5}$. For $p=5$, the $L^2$-distance $\| u_p-u_{\text{opt}}\|_{L^2(0,T)}$ is approximately 10 times smaller than for $p=2$. For $\beta= 10^{-4}$ and $\beta= 5 \cdot 10^{-5}$, a significant reduction of the costs can be observed as $k$ increases.
For $\beta= 10^{-5}$, the situation is more complex. Convergence of the controls $u_k$ as $k$ increases is not achieved.
The values of the costs decrease with $k$ except for $k=4$. In case the associated closed loop system associated with still converges to $0$ but some strong
oscillations of the control render a large value of $J(u_4)$.
% It is not clear if these variations of
%the control are intrinsic to the feedback law or if they are due to
%approximation errors occurring when simulating the closed-loop system, which is
%highly non-linear.
%The control $u_5$ appears to be again very efficient, in comparison with $u_2$.

\paragraph{Test case 5: initial condition with support in the second potential-well}

Figure \ref{fig:2d_gaussian} and Table \ref{table:2d_gaussian} show the results
obtained for an initial condition located in the second potential well (around
$x_B$). A large proportion of the distribution must be transported to $x_A$,
along the first coordinate axis and from the negative values to the positive
ones.
Therefore, one can expect that the first coordinate of the control takes
negative values and that the second coordinate has a smaller amplitude than the first one.
This initial condition is much further from the stationary distribution than  the previous one. Consequently the cost of the uncontrolled system is much larger: $J(u=0)= 9.04$ and $\| \rho_0 - \rho_\infty
\|_{L^2(\Omega)} = 0.63$. As a consequence, the feedback laws are only efficient
for larger values of $\beta$ than those considered previously. Convergence can
be observed for values of $\beta$ larger than $5 \cdot 10^{-3}$. The reduction
factor of the $L^2$-distance is smaller than for the test case 4, but still a significant
reduction of the cost is noted for $\beta= 10^{-2}$ and $\beta= 5 \cdot
10^{-3}$. For $\beta= 10^{-3}$,
convergence with respect to $k$ cannot be achieved. The closed-loop system associated with
$\mathbf{u}_4$ quickly converges to a non-trivial stationary point. The
closed-loop system associated with $\mathbf{u}_5$ generates a control which has
strong oscillations  along time and eventually converges to a non-trivial
stationary point.

\section{Conclusion}

A numerical method for computing polynomial feedback laws for an
infinite-dimensional optimal control problem with infinite-time horizon has been
proposed. It consists in particular in reducing the state equation in order to
attenuate the curse of dimensionality, which prevents a direct resolution of the
involved Lyapunov equations. The applicability of the method has been
demonstrated with an optimal control problem of the Fokker-Planck equations in
dimensions 1 and 2. The effect of model reduction on the feedback laws has been
numerically analysed and the relevance of the reduction approach has been shown.
Good convergence results for high-order feedback laws have been obtained in many
situations for which the initial condition was close enough to the equilibrium
or for which the value of the cost parameter $\beta$ was not too small. The influence of
the initial condition and the cost parameter $\beta$ on the success of the
method has been investigated in a systematic manner.

Further research will focus on the design of polynomial feedback laws for
infinite-dimensional systems with a more complicated structure. At a numerical
level, the use of low-rank tensors formats could be investigated to facilitate
the numerical resolution of the Lyapunov equations and the simulation of
closed-loop systems. It may also be of interest to design a heuristic mechanism
which selects an appropriate order for the feedback law, in order to avoid
convergence to a non-trivial stationary point and to allow a practical
implementation. At a theoretical level, the computation of an error estimate for
the efficiency of controls generated by reduced feedback laws could also be a
topic for future work.

\section*{Acknowledgements}

This work was partly supported by the ERC advanced grant 668998 (OCLOC) under the EU's H2020 research program.

%%%%%%%%%

\end{document}